%% file: Analysis.tex
\documentclass[11,letterpaper,leqno]{article}
\usepackage[letterpaper, hmargin=2cm, vmargin={2cm, 2cm}]{geometry}
\usepackage{myDefs}

\begin{document}

\title{Analysis of Boundary Behaviour of Quasidisks and Jordan Repellers}
\author{Ilia Binder, Adi Gl\"ucksam}
\date{}
\maketitle
\begin{abstract}
We investigate the fine properties of harmonic measure and boundary rotation. By focusing on quasidisks and, in particular, on connected Jordan Repellers arising from conformal expanding dynamical systems, we explore those using the deep interplay between geometric function theory and dynamical systems as a unified framework.\\
    While some of the results presented here are regarded as 'folklore' among experts, they lack rigorous proofs in the existing literature. We fill this gap by providing a comprehensive, referable treatment using a novel approach that also expands existing results.
\end{abstract}

\section{Introduction}\label{sec:intro}
\subfile{sections/introduction}

\section{Background, Definitions, and Results}\label{sec:defs}
\subfile{sections/definitions}

\section{Preliminaries: Conformal Maps, Harmonic Measure, and Rotations}\label{sec:preliminaries}
\subfile{sections/preliminaries}

\section{Derivative vs. Harmonic Measure and Rotation}\label{sec:relations}
\subfile{sections/relations}

\section{Minkowski Crosscut Spectrum}\label{sec:curve_spec}
\subfile{sections/spectra}

\section{Multifractal Analysis and Thermodynamic Formalism}\label{sec:repellers}
\subfile{sections/repellers}

\bibliographystyle{plain}


\end{document}

%% file: sections/introduction.tex
It is by now a classical result that the boundary behavior of conformal maps and the dimension properties of harmonic measure of simply connected domains are intimately related. This connection is profoundly illustrated by Makarov's celebrated theorem on the dimension of the harmonic measure for simply connected planar domains, \cite{MakThm}. The proofs of this theorem typically establish that the derivative of a conformal map grows sub-exponentially at almost every point of the boundary. Further insight into this relationship is provided by the celebrated paper, \cite{CM}, which studies negative moments of the derivatives of conformal maps using estimates on the number of disjoint disks with large harmonic measure. Less widely utilized is a similar relationship between the local rotation of a simply connected domain's boundary and the growth of the argument of its corresponding conformal map's derivative. These relationships are best formalized using \emph{Multifractal analysis}.

Roughly speaking, the \emph{distortion spectrum}, denoted $d(a,b)$, describes the dimension of the set of points on the circle with a given growth-rate of the derivative and its argument. The entropy function for the distortion spectrum, is the \emph{integral means spectrum}, which conceptually describes the \emph{mean} growth-rate of the derivative of a conformal map raised to a given power $z$ as we approach the boundary. We allow the power $z$ to be a \emph{complex} number. While in the context of real exponents, this object is classical and has been extensively studied (e.g., \cite{makbook} and \cite{pombook}), however, in the context of complex exponents, there are still many important open questions, such as Brennan's Conjecture. Values of the integral mixed spectrum for real exponents reflect the behavior of the harmonic measure, while purely imaginary exponents correspond to the behavior of the rotation. Most of the classical results can be easily extended to this version of the integral means spectrum; nevertheless, it gives rise to many new and intriguing questions, see, e.g., \cite{Baranov2008}. We informally refer to those spectra as 'conformal'.

On the potential side, the counterpart of the distortion spectrum is the \emph{dimension spectrum}. Similarly to the distortion spectrum, intuitively,  the \emph{dimension spectrum} of a measure, denoted $f(\alpha)$, aims to describe the dimension of the set of points where the measure exhibits a local scaling behavior corresponding to dimension $\alpha$. The counterpart of the integral mean spectrum, i.e., the entropy function for the dimension spectrum is the packing spectrum. For a detailed discussion of Multifractal analysis, we refer the reader to \cite{fal}. 

In this paper, we study these spectra of harmonic measures, and explore similar notions for the rotation of the local boundary, which necessitates the introduction of a two-parameter version of the dimension spectrum, $f(\alpha, \gamma)$ (see Sections \ref{sec:defs} for precise definition). Multifractal analysis of harmonic measure was investigated in the one-parameter case in Makarov's celebrated paper, \cite{makbook}. We informally refer to those multifractal spectra as 'geometric'.

 The previously mentioned entropy relations manifests as a Legendre-type relation between the integral means spectrum and the distortion spectrum (specifically, a version of the  Minkowski distortion spectrum) and between the packing spectrum and the dimension spectrum. These relations will not be explored in this note. We refer to \cite{dis, makbook} for further details. 

The relationship between the 'geometric' and 'conformal' spectra is substantially more intricate than expected. From a 'naive' understanding of conformal maps, one expects the following relationship:
\begin{equation}\label{eq:distdim}
d_\Om(a,b)= (1-a)f_\Om\left (\frac1{1-a}, \frac{-b}{1-a}\right).\end{equation}
Nevertheless, somewhat surprisingly, the examples presented in \cite{makbook} and \cite{AdiBin2} illustrate that this relationship does not hold, even for John domains. More precisely, these examples demonstrate that  either spectrum can be strictly dominant for all domains.

On the other hand, in the specific case of quasidisks, domains bounded by quasicircles, Relation \eqref{eq:distdim} holds true, as outlined in Theorem \ref{thm:fine=distortion_quasidisks}. While this theorem is widely recognized as a 'folklore result' among experts in the field, we are not aware of any complete published proof. One of the primary objectives of our paper is to fill this notable gap in the literature, an undertaking that proved to be more complex than expected, as elaborated in Section \ref{sec:relations}.

In Section \ref{sec:curve_spec}, we introduce an innovative notion that can be thought of as the 'crosscut' version of the Minkowski dimension spectrum. This new version extends the conventional framework and offers new insights, bridging, in some sense, the 'geometric' and the 'conformal' worlds. An essential property of the crosscut spectrum is that it majorates the usual \emph{distortion} spectrum for arbitrary domains and aligns with it for quasidisks. We utilize this version of the spectrum in \cite{AdiBin2} to produce counterexamples and to generate an upper bound for the universal distortion spectrum.

Multifractal analysis offers a particularly elegant framework for a subclass of quasidisks: Jordan repellers (see \cite{makbook} and Section \ref{subsec:repellers} for the precise definition). Jordan repellers can also be described using expanding dynamical systems - shifts over finite alphabets. We introduce parallel definitions of spectra using the framework of words over finite alphabets and demonstrate their equivalence to the original ones. The fact that Jordan repellers are invariant under an expanding dynamics allows us to utilize the power of Thermodynamic Formalism. Using this well-established technique, we show that for repellers, all versions of the spectra - dimension, and distortion as well as Hausdorff, and Minkowski - coincide, and the dimension and distortion spectra are related via Legendre transforms to each other. We provide proofs of these assertions in the case of complex exponents (for the conformal spectra) and the two-parameter geometric spectra, extending the results of \cite{makbook}. While all of these results can also be considered 'folklore', we need a stronger version of these relations in \cite{AdiBin2}. These stronger relations come from our extension of the classical Carleson estimate on multiplicativity of harmonic measure, Lemmas \ref{lem:Carleson} and \ref{lem:Carleson_rotations}. 

The other objective of this note is to prepare the necessary auxiliary results for \cite{AdiBin2}. In \cite{AdiBin2}, we show, among other things, that unlike the standard spectra, measured for a specific set, the universal conformal and geometric spectra mentioned above \emph{do} satisfy Relation \eqref{eq:distdim}. This is yet another 'folklore result'. We establish it by showing that the universal bounds can be obtained by looking solely at \emph{Jordan repellers}.

Section \ref{sec:repellers} contains 
one important technical result; The introduction of a \emph{'propagating property'}, which, roughly speaking, is a property whose density does not decrease as the scale shrinks. We prove that the spectra of repellers propagate, enabling controlled approximations of the spectra by using spectra at a specific scale - finer scales yield better approximations. We expect the concept of propagating property to be fruitful in solving other problems as well.

Many of the claims proven in this paper are intuitive and fundamental to the field; however, we were unable to find references or published proofs in standard literature. We have therefore included full proofs both for completeness and to provide a convenient reference for future use.

%% file: sections/definitions.tex
\subsection{Setup and Fundamental Spectra}
Let $\Omega$ be a bounded simply connected domain, fix $z_0\in\Omega$, and let $\phi:\D\to\Omega$ be a conformal map sending the origin to $z_0$. Let $\arg$ denote the branch of the argument of $\phi'$ satisfying $\arg(\phi'(0))\in(-\pi,\pi]$.

We first discuss the 'conformal spectrum':
\begin{defn}
For every $ a, b\in\R$, the \emph{Minkowski Distortion Mixed Spectrum} of $\Omega$ is defined by
$$
d^{\sigma\sigma'}_\Omega(a,b):=\limsup_{\eta\searrow 0}\;\limsup_{r\nearrow1^-}\frac{\log\bb{\lambda_1\bb{L^{\sigma\sigma'}_{a,b,\eta}(r)}}}{\log\bb{\frac1{1-r}}}+1
,\quad \sigma,\sigma'\in\bset{+,-}$$
where $\zeta\in L^{\sigma\sigma'}_{a,b,\eta}(r)\subset\partial\D$ if 
\begin{center}
\begin{minipage}[b]{.5\textwidth}
\vspace{-\baselineskip}
\begin{equation*}
\tag{d1}\label{eq:d_abs}
\abs{\phi'(r\zeta)}\begin{cases}>\bb{1-r}^{-a+\eta}&, \sigma=+\\
							<\bb{1-r}^{-a-\eta}&, \sigma=-
				\end{cases}.
\end{equation*}
\end{minipage}%
\begin{minipage}[b]{.5\textwidth}
\vspace{-\baselineskip}
\begin{equation*}
\tag{d2} \arg\bb{\phi'(r\zeta)}\begin{cases}>\bb{1-r}^{-b+\eta}&, \sigma'=+\\
							<\bb{1-r}^{-b-\eta}&, \sigma'=-
				\end{cases}. \label{eq:d_rot}
\end{equation*}
\end{minipage}
\end{center}
The following version of the spectrum, which is essentially the minimum of the four quantities defined above, and is also extremely useful.
$$
d_\Omega(a,b):=\limsup_{\eta\searrow0}\;\limsup_{r\nearrow1^-}\frac{\log\bb{\lambda_1\bb{L_{a,b,\eta}(r)}}}{\log\bb{\frac1{1-r}}}+1,
$$
where $\zeta\in L_{a,b,\eta}(r)\subset\partial\D$ if
$$
\bb{1-r}^{-a+\eta}<\abs{\phi'(r\zeta)}<\bb{1-r}^{-a-\eta}\quad\quad;\quad\quad  \bb{1-r}^{-b+\eta}<\arg\bb{\phi'(r\zeta)}<\bb{1-r}^{-b-\eta}.
$$
\end{defn}
Given a function $g:\R^2\rightarrow\R$ which is increasing (decreasing) in each variable, we say a point $(a,b)$ is a \emph{growth point} if for all $a$ and $\eta>0$, $g(a+\eta,b)>g(a-\eta,b)$ ($g(a+\eta,b)<g(a-\eta,b)$) and for all $b$ and $\eta>0$, $g(a,b-\eta)<g(a,b+\eta)$ ($g(a,b+\eta)<g(a,b-\eta)$). By definition,
$$
d_\Omega(a,b)\le\min\bset{d_\Omega^{++}(a,b),d_\Omega^{+-}(a,b),d_\Omega^{--}(a,b),d_\Omega^{-+}(a,b)},
$$
and if $(a,b)$ is a growth point for $d_\Omega^{\sigma,\sigma'}( a,b)$ then $d_\Omega(\alpha,\gamma)=d_\Omega^{\sigma,\sigma'}(\alpha,\gamma)$.

Note that the spectrum does not depend on the choice of $z_0$ or $\phi$. Moreover, by Makarov's triple logarithm theorem,  $d_\Omega(0,0):=1$, since almost everywhere $\log\phi'$ grows like $o\left(\frac1{1-r}\right)$.  Roughly speaking, the Minkowski Distortion Spectrum measures the size of the set of 'atypical' points on the circle with fast growth or decay of $\abs{\phi'}$ and $\arg\bb{\phi'}$. We note that this definition appears to differ from the one given for the non-rotational case in \cite{makbook}. The reason for the more complicated definition lies in the fact that for fixed $a\neq0$, the function $d_\Omega(a,b)$ might not achieve its maximum at $b=0$.

Now, let us turn to the 'geometric spectra'. To discuss the \emph{geometric} behavior of the boundary, we first need to define rotation. We define the \emph{rotation} of a domain, $\Omega$, centered at a boundary point $z$ at scale $\delta$ by
$$
\rot(z,\delta):=\exp\bb{\underset{y\in\partial\Omega_{z,\delta}\cap \partial B(z,\delta)}\inf \arg_{[z]} (y-z)},
$$
where $\Omega_{z,\delta}$ is the connected component of the set $\bset{y\in\Omega,\;\abs{y-z}>\delta}$ containing $z_0$, and the argument, $\arg_{[z]}$, is a branch of the argument satisfying $\arg_{[z]}(z_0-z)\in(-\pi,\pi]$.

It will sometimes be useful to discuss the rotation of a disk, $B$. We will abuse the notation of rotation, $\rot$, and define the rotation of a disk  $B=B(z,\delta)$ by $\rot(B)=\rot(B(z,\delta)):=\rot(z,\delta)$, even when its centre is not necessarily on $\partial\Omega$ as long as $\frac{diam(\partial\Omega\cap B)}{diam(B)}\geq\frac14$.

\begin{defn}
For every $\alpha,\gamma\in\R$ the  \emph{Minkowski Dimension Mixed Spectrum} of $\Omega$ is defined by
$$
f_\Omega^{\sigma\sigma'}(\alpha,\gamma)=\limdec \eta 0 \limitsupdec \delta 0 \frac{\log N^{\sigma\sigma'}(\delta,\alpha,\gamma,\eta)}{\log\bb{\frac1\delta}},\quad \sigma,\sigma'\in\bset{+,-}
$$
where $N^{\sigma\sigma'}(\delta,\alpha,\gamma,\eta)$ is the maximal number of pairwise disjoint disks, $\bset{B(z_j,\delta)}$, centered at $\partial\Omega$, for which
\begin{center}
\begin{minipage}[b]{.5\textwidth}
\vspace{-\baselineskip}
\begin{equation*}
\tag{f1}\label{eq:N_measure}
\omega(B(z_j,\delta))\begin{cases} >\delta^{\alpha+\eta}&, \sigma=+\\
								<\delta^{\alpha-\eta}&, \sigma=-
					\end{cases}.
\end{equation*}
\end{minipage}%
\begin{minipage}[b]{.5\textwidth}
\vspace{-\baselineskip}
\begin{equation*}
\tag{f2} \label{eq:N_rotation}
\rot(z_j,\delta)\begin{cases}>\delta^{\gamma+\eta}&, \sigma'=+\\
							<\delta^{\gamma-\eta}&, \sigma'=-
				\end{cases}.
\end{equation*}
\end{minipage}
\end{center}

As we have done above for the distortion spectrum, we also define 
$$
f_\Omega(\alpha,\gamma):=\limdec \eta 0 \limitsupdec \delta 0 \frac{\log N(\delta,\alpha,\gamma,\eta)}{\log\bb{\frac1\delta}},
$$
where $N(\delta,\alpha,\gamma,\eta)$ is the maximal number of pairwise disjoint disks $\bset{B(z_j,\delta)}$ centered at $\partial\Omega$ satisfying that
$$
\delta^{\alpha+\eta}< \omega(B(z_j,\delta))< \delta^{\alpha-\eta}\quad\quad;\quad\quad \delta^{\gamma+\eta}<\rot(z_j,\delta)<\delta^{\gamma-\eta}.
$$
\end{defn}
And, as before,
$$
f_\Omega(\alpha,\gamma)\le\min\bset{f_\Omega^{++}(\alpha,\gamma),f_\Omega^{+-}(\alpha,\gamma),f_\Omega^{--}(\alpha,\gamma),f_\Omega^{-+}(\alpha,\gamma)},
$$
and if $(\alpha,\gamma)$ is a growth point for $f_\Omega^{\sigma,\sigma'}(\alpha,\gamma)$ then $f_\Omega(\alpha,\gamma)=f_\Omega^{\sigma,\sigma'}(\alpha,\gamma)$.

\begin{rmk}\label{rmk:adjoint_set}
Let $\phi:\D\rightarrow\Omega$ be a conformal map and denote by $\tilde\phi(z):=\overline{\phi\bb{\overline z}}$. Since $\tilde \phi'(z)=\overline{\phi'\bb{\overline z}}$, $\tilde\phi$ is also conformal, and its image, $\tilde\Omega$, satisfies
$$
\tilde\Omega:=\bset{\tilde\phi(z), z\in\D}=\bset{\overline{\phi\bb{\overline z}}, z\in\D}=\bset{\overline{\phi\bb{z}}, z\in\D}=\overline\Omega.
$$
In particular, $\tilde{\tilde\Omega}=\Omega$. We conclude that for every $a,b,\alpha,\gamma$ we have
$$
d_\Omega^{\pm+}(a,b)=d_{\tilde \Omega}^{\pm-}(a,-b)\quad, \quad f_\Omega^{\pm+}(\alpha,\gamma)=f_{\tilde\Omega}^{\pm-}(\alpha,-\gamma).
$$
We will therefore only consider $b\ge 0$ or $\gamma\ge 0$, as the general case can be obtained using these cases and the equalities above.
\end{rmk}

\subsection{Relations Among the Fundamental Spectra}
The first significant result in this note is a theorem which states that if the harmonic measure is doubling (or equivalently, the domain is a quasidisk), then the Minkowski Dimension Mixed Spectrum is directly related to the Minkowski Distortion Mixed Spectrum (with the correct parameter).

\begin{thm}\label{thm:fine=distortion_quasidisks}
Let $\Omega\subset\C$ be a quasidisk. Then
\begin{enumerate}
\item\label{item:d>f} $\quad\quad\quad\quad\quad\quad\quad\quad\quad d_\Omega(a,b)=\bb{1-a}  f_\Omega\bb{\frac1{1-a},\frac{-b}{1-a}}.$
\item\label{item:d<f} $\quad\quad\quad\quad\quad\quad\quad\quad\quad \min_{\sigma,\ \sigma'\in\bset{-,+}}\bset{d^{\sigma\sigma'}_\Omega(a,b)}=\bb{1-a}  \min_{\sigma,\ \sigma'\in\bset{-,+}}\bset{f^{\sigma\sigma'}_\Omega\bb{\frac{1}{1-a},\frac{-b}{1-a}}}.$
\end{enumerate}
\end{thm}

For approximately thirty years, experts in the field believed that the inequality 
\begin{equation}\label{eq:d>f}
d_\Omega(a,b)\geq\bb{1-a}  f_\Omega\bb{\frac1{1-a},\frac{-b}{1-a}}
\end{equation}
holds without the assumption that $\Omega$ is a quasidisk. While attempting to write down a rigorous proof, the authors encountered a problem - it is possible that every curve in a disk $B$ either carries a large portion of the harmonic measure or has a large enough diameter, but not both. Extending the curve would either increase the harmonic measure too much or make the curve too long. Such disks should be counted on a different, smaller scales. 

As established in this paper, the issue does not occur for quasidisks. However, the authors constructed a simply-connected John domain for which the Minkowski Distortion Spectrum is smaller than the Minkowski Dimension Spectrum, i.e., inequality \eqref{eq:d>f}, and consequently Theorem \ref{thm:fine=distortion_quasidisks}, do not hold (see \cite{AdiBin2}). Note that the very existence of this counterexample necessitated a new approach for the proof of the famous fractal approximation, which is also covered in \cite{AdiBin2}.

In order to prove Theorem \ref{thm:fine=distortion_quasidisks}, we rely on two 'folklore results' that are of independent interest. To state those results, we require notation, some of which is standard in the field and some is new.

\paragraph*{A Word About Notation}
\begin{itemize}
\item  Throughout the paper, we write $a\sim b$ if there are constants $0<c_1<c_2<\infty$ so that $c_1\le\frac ab\le c_2.$
We write $a\dsim W b$ if the constants depend on the object $W$.\\We extend this notation to accommodate inequalities as well; we write $a\lesssim b$ if there exists a constant $c>0$ so that $a\le c\cdot b$ and $a\dless W b$ if this constant depends on the object $W$.
\item Given a disk $B=B(z,\delta)$ the set $L\cdot B$ denotes the disk concentric with $B$ of radius $L\cdot\delta$, i.e., $L\cdot B =B(z,L\cdot\delta)$.
\item Following \cite[p.6]{pombook} we will denote by $\rho(z,w)=\tanh\;\lambda(z,w)=\abs{\frac{z-w}{1-\overline z\cdot w}}$, there $\lambda$ is the hyperbolic metric.
\end{itemize}
\begin{defn}\label{def:represent}
Let $I\subset\partial \D$ be an arc satisfying $\lambda_1(I)<\frac12$. Unless otherwise stated, we denote by $\zeta_I$ the centre of the arc, $I$, and by  $z_I=\zeta_I\bb{1-\lambda_1(I)}$. We refer to $z_I$ as \emph{the point representing  $I$}. Note that this representation may not be a 'good representation', but, as we will see, for quasidisks, it is always a 'good representation'.

Given a disk, $B$, by Carleson's lemma, \cite[Lem 2.5 p. 277]{gm}, there exists a continuum $\beta_B=\beta\subset\partial\Omega\cap 2B$, satisfying that $\omega(\beta)\in\bb{\frac{\omega(B)}{\log^2\bb{\frac1{\omega(B)}}},\omega(B)}$. We denote by $z_B=\bb{1-\omega(B)}\cdot\zeta_B$ where $\zeta_B$ is the center of the arc $\phi\inv(\beta)$. We refer to $\beta_B$ as \emph{the arc representing $B$} and to $z_B$ as \emph{the point representing  the disk $B$}. Note that as in the case of an arc, this representation may not be a 'good representation', but for quasidisks it is always a 'good representation'.
\end{defn}

Lastly, for every Jordan curve, $C\subset\C$, and every $w_1,w_2$ we denote by $C_{w_1,w_2}\subset C$ the connected component of $C\setminus\bset{w_1,w_2}$ of smaller diameter. Given a quasidisk, $\Omega$, we denote by $K_\Omega$ the dilatation of $\Omega$ (for the precise definition see Definition \ref{def:QD}).

\begin{thm}[The Relation Theorem]\label{thm:relations}
Let $\Omega\subset\C$ be a quasidisk.  We consider the properties
\begin{center}
\begin{minipage}[b]{.5\textwidth}
\vspace{-\baselineskip}
\begin{equation*}
\tag{HM} \label{eq:derivative_measure} \frac{\delta}{\omega(B)\cdot \log^2\bb{\frac1{\omega(B)}}}\dless{K_\Omega}\abs{\phi'\bb{z}} \dless{K_\Omega} \frac{\delta}{\omega(B)}
\end{equation*}
\end{minipage}\hfill
\begin{minipage}[b]{.4\textwidth}
\vspace{-\baselineskip}
\begin{equation*}
\tag{R}  \label{eq:derivative_rotation} \abs{\phi'^{-i}\bb{z}}\dsim{K_\Omega} \rot(B).
\end{equation*}
\end{minipage}
\end{center}
Then-
\begin{enumerate}
\item 
For every disk, $B$, of radius $\delta$ centered at $\partial\Omega$, the point $z_B$, which represents the disk $B$, satisfies \eqref{eq:derivative_measure} and \eqref{eq:derivative_rotation}.  If $2B\cap 2B'=\emptyset$ then $\rho(z_B,z_{B'})\dmore{K_\Omega}1$.
\item
For every point $z\in\D\setminus\{0\}$ let $A_z:=\bset{\xi\in\partial\D, \abs{\xi-\frac z{\abs z}}<1-\abs z}$. There exists a disk $B_z$ of radius $\delta\overset{K_\Omega}\sim\bb{1-\abs z}\cdot\abs{\phi'(z)}$ centered at $\partial\Omega$, containing $\phi(A_z)$ and satisfying \eqref{eq:derivative_measure} and \eqref{eq:derivative_rotation}. There exists $\alpha(K_\Omega)$ so that if $\rho(z,z')\ge\alpha(K_\Omega)$ then $B_z\cap B_{z'}=\emptyset$.
\end{enumerate}
\end{thm}

\begin{defn}
A \emph{crosscut} is a continuous injective map $\gamma:[0,1]\rightarrow\overline\Omega$, satisfying that $\gamma(0),\gamma(1)\in\partial\Omega$ while $\gamma(0,1)\subset\Omega$. The \emph{support of a crosscut}, $\gamma$, is the set $\nu_\gamma=\partial\bb{\Omega_\gamma}\setminus \gamma\subset\partial\Omega$, where $\Omega_\gamma\subset\Omega$ is the connected component of $\Omega\setminus \gamma$ which does not contain $z_0$.\\
Given an arc $I\subset\partial\D$ we denote by $C_I\subset\Omega$ a crosscut satisfying
\begin{enumerate}
\item \label{itm:large_cut}$\phi\inv(C_I)\subset\D$ is a crosscut supported on some set $J \subset\partial\D$, with $I\subseteq J\subset I\bb{1+10^{-3}}$.
\item For any other crosscut, $C$,  satisfying  \ref{itm:large_cut}, $\diam(C_I)\le 2\diam(C)$.
\end{enumerate}
We call $C_I$ \emph{a crosscut of $I$}. We will sometimes abuse the notation $\nu_I$ to denote the support of a crosscut $C_I$.
\end{defn}
Note that $C_I$ is not uniquely defined, and in fact, it cannot be uniquely chosen in some cases, e.g., when one of the end-points of $I$ is not degenerate.
\begin{rmk}
Why every arc $I$ has a crosscut, $C_I$?\\
While if $C$ is a crosscut in $\Omega$, then $\phi\inv(C)$ is a crosscut in $\D$ (see \cite[Proposition 2.14]{pombook}), it could be the case that $C$ is a crosscut in $\D$ but $\phi(C)$ is no longer a crosscut in $\Omega$. However, if $C\subset\D$ is a crosscut whose end-points have radial limits, then $\phi(C)\subset\Omega$ is a crosscut as well (see e.g. \cite[Theorem 2.17(1)]{pombook}). Since the radial limits of $\phi$ exists at almost every point on $\partial\D$ (see  \cite[Theorem 1.7]{pombook}), for every arc $I\subset\partial\D$ a crosscut of $I$ exists.\\
Note that $\diam(\nu_{C_I})$ might be much bigger than $\diam(C_I)$, e.g., if $\nu_{C_I}$ is a cusp. However, $\diam(\nu_{C_I})$ always dwarfs $\frac12\diam(C_I)$.
\end{rmk}

\begin{defn}\label{defn:rotation_curve}
We define the \emph{Rotation} of a crosscut, $C$, by
$$
\rot(C)=\underset{\xi\in \nu_C}\inf \rot(\xi,\diam(C)).
$$
\end{defn}
Note that since the end-points of $C$ are on $\partial\Omega$, then for every $\xi\in C$, $\dist(\xi,\partial\Omega)\le\frac{\diam(C)}2$ implying that
$$
\frac{\diam(\partial\Omega\cap B(\xi,\diam(C)))}{\diam(B(\xi,\diam(C)))}\ge\frac{\bb{2-\frac12}\diam(C)}{2\diam(C)}=\frac34>\frac14.
$$
We see that the rotation $\rot(\xi,\diam(C))$ is well defined for every $\xi\in C$. In addition, we note that for Jordan domains, the support of a crosscut is always a curve.

The proof of The Relation Theorem, Theorem \ref{thm:relations}, relies on  the following Lemma, which holds for general domains. 

\begin{lem}[The Main Lemma]\label{lem:main}
Let $\Omega\subset\C$ be a simply connected domain and let $\phi:\D\rightarrow\Omega$ be a Riemann map sending $0$ to $z_0$. Then
\begin{itemize}
\item[]
\vspace{-1cm}
\begin{equation*}
\hspace{-6cm}\tag{HM*}\label{eq:arc_measure}
\abs{\phi'\bb{z_I}} \sim\frac{\dist(\phi(z_I),\partial\Omega)}{\lambda_1(I)}\sim \frac{\diam(C_I)}{\lambda_1(I)}\lesssim \frac{\diam(\nu_I)}{\lambda_1(I)}.
\end{equation*}
\item []
\vspace{-1cm}
\begin{equation*}
\hspace{-10.5cm}\tag{R*} \label{eq:arc_rotation} \lim_{r\to1^-}\frac{\log\abs{\phi'^{-i}\bb{z_I}}}{\log\rot(C_I)}=1.
\end{equation*}
\item [] 
\vspace{-1cm}
 \begin{equation*}
\hspace{-9.5cm} \tag{D*}\label{eq:disjoint_dsks} \rho(z_B,z_{B'})\sim \frac{\omega\bb{C_{I_{B,B'}}}}{\min\bset{\omega(B),\omega(B')}},
\end{equation*}
for every disjoint disks $B,B'$ centered at $\partial\Omega$, and $I_{B,B'}=\sbb{\frac{z_B}{\abs{z_B}},\frac{z_{B'}}{\abs{z_{B'}}}}\subset\partial\D.$
\end{itemize}
If $\Omega$ is a quasidisk, those properties become stronger. Namely, if $K_\Omega$ is the dilatation of $\Omega$, then
\begin{itemize}[]
\addtolength{\itemindent}{1.2cm}
\item[(HM*)]  $\abs{\phi'\bb{z_I}}\dsim{K_\Omega} \frac{\diam(\nu_I)}{\lambda_1(I)}$.
\item [(R*) ~ ]  $ \abs{\phi'^{-i}\bb{z_I}}\dsim{K_\Omega} \rot(C_I)$.
\item [(D*) ~ ]  If $2B\cap 2B'=\emptyset$, then $\rho(z_B,z_{B'})\overset{K_\Omega}\gtrsim1$.
\end{itemize}
\end{lem}

\subsection{Corsscut Spectrum}
Attempting to relate the 'geometric' and the 'conformal' spectra, we define the \emph{Minkowski Crosscut Spectrum}.
\begin{defn}\label{def:curves_and_d}
For every $r\in\bb{0,1}$ and $a,b>0$ we denote by $\Gamma\bb{a,b,r}$ a maximal collection of crosscuts, $C$, with disjoint supports satisfying
\begin{enumerate}
\item $\abs{\frac{\omega(\nu_C)}{1-r}-1}\le\frac{1-r}{\log\log\log\bb{\frac1{1-r}}}$.
\item $\diam(\nu_C)\ge\bb{1-r}^{1-a}$.
\item $\rot(C)\ge\bb{1-r}^{-b}$.
\end{enumerate}
We define the \emph{Minkowski Crosscut Spectrum} by
$$
d_\Omega^{cc}(a,b)=\underset{a'\searrow a\atop {b'\searrow b}}\limsup\;\underset{r\nearrow 1}\limsup\;\frac{\log\bb{\# \Gamma\bb{a',b',r}}}{\log\bb{\frac1{1-r}}}.
$$
\end{defn}
We relate this new notion of spectrum to the fundamental ones discussed in previous sub-sections:
\begin{lem}\label{lem:curves_bnd_d}
Let $a>0, b\in\R$.
\begin{enumerate}
\item  \label{item:ineq_ds} For every simply connected domain, $\Omega$, $d_\Omega(a,b)\le d_\Omega^{++}(a,b)\le d_\Omega^{cc}(a,b)$.
\item\label{item:quasidisk_curve} For every quasidisk, $\Omega$, $d_\Omega(a,b)= d_\Omega^{cc}(a,b)$.
\end{enumerate}
\end{lem}

\begin{lem}\label{lem:neg_a}
If $a<0$, then for every simply connected domain, $\Omega$,
$$
d_\Omega(a,b)\le d_\Omega^{++}(a,b)\le \limitinf\eta{0^+} f_\Omega^{+\pm}\bb{\frac{1}{1-a}+\eta,\frac{-b}{1-a}+\eta}.
$$
\end{lem}

\subsection{Jordan Repellers}\label{subsec:repellers}
The last section of this note deals with the special case of \emph{Jordan Repellers}.
\begin{defn}
Let $\Omega$ be a simply connected domain, and $\partial$ be a sub-arc of its boundary, $\partial\Omega$. We call $\partial$ a \emph{Jordan Repeller} if there exists a partition of $\partial$ into a finite number of non-intersecting sub-arcs, denoted $\partial_1,\partial_2,\cdots,\partial_N$ (a Markov partition), and a piecewise univalent (though not necessary injective) map $F$, such that for every sub-arc in the partition, $\partial_j$, there exists a neighborhood $U_j$ containing $\partial_j$ and a univalent map $F_j:\, U_j\rightarrow\C$ satisfying:
\begin{enumerate}[label=(\arabic*)]
\item $F|_{\partial_j}=F_j$.
\item (Geometry invariance) $F_j\bb{\Omega\cap U_j}\subset\Omega$. (This property preserves the geometry of the set).
\item (Boundary invariance) $F_j\bb{\partial\Omega\cap U_j}\subset\partial\Omega$.
\item (Markov property) The image of each $\partial_j$ is a finite union arcs $\partial_k$, i.e., $F(\partial_j)=\underset{k\in\mathcal A_j}\bigcup\partial_k$, where $\mathcal A_j\subseteq\bset{1,2,\cdots,N}$.
\item (Expanding) There exists $n_0$ so that $\underset{z\in\partial}\inf\;\abs{F^{(n_0)} (z)} >1$, where $F^{(n_0)}$ denotes the $n_0$'th derivative of $F$.
\item (Mixing) For any open set $W$ satisfying $W\cap\partial\neq\emptyset$, there exists $n$ such that $F^{\circ n}\bb{W\cap\partial}=\partial$, i.e., the $n$'th iteration of $F$ maps $W\cap\partial$ to cover $\partial$.
\end{enumerate}
\end{defn}


The main result of Section \ref{sec:repellers} is a refinement and an extension of Carleson's estimate to explore the multiplicativity properties of rotations and harmonic measure of Jordan Repellers (Lemmas \ref{lem:Carleson} and \ref{lem:Carleson_rotations}).

Finally, we present the notion of 'propagation':
\begin{defn}
A function $\varphi:(0,1)^2\rightarrow\R$ \emph{propagates} if there exists a constant $C>1$, so that for every $n\in\N$
$$
\frac{\log \varphi\bb{\delta^{n},C\cdot \eta}}{\log\bb{\frac1{\delta^{n}}}}\ge \frac{\log \varphi(\delta,\eta)}{\log\bb{\frac1\delta}}
$$
for every $\delta$ small enough (which may depend on $\eta$).
\end{defn}
In section \ref{sec:repellers}, we show that the spectra of Jordan Repellers propagate. This crucial property allows us to lay the foundations to approximate the universal spectra via Repellers.

%% file: sections/preliminaries.tex
In this section, we provide the necessary background on conformal maps, harmonic measures, and rotations. We include full proofs of several key results here, as published documentation for intuitive or 'folklore claims' which is absent from the standard literature. This ensures completeness and provides a convenient reference.

The first subsection focuses on quasidisks and contains results concerning the interrelations between conformal maps and harmonic measures in this special setup. The second subsection details the results concerning rotations.

\subsubsection{On Two Notions of Doubling Measures on Quasicircles}
\begin{defn}\label{def:QD}
Following \cite[p. 234-235]{gm}, we say that $\Gamma$ is a \emph{quasicircle} if there exists a constant $K$ such that if $C_{w_1,w_2}\subset\Gamma$ is the curve of smallest diameter connecting $w_1$ and $w_2$, then for every $w\in C_{w_1,w_2}$
$$
\frac{\abs{w_1-w}+\abs{w_2-w}}{\abs{w_1-w_2}}<K.
$$
A \emph{quasidisk} is a domain $\Omega$, bounded by a quasicircle. We will call $K$ the \emph{dilatation} of $\Omega$. Although $K$ is not precisely the dilatation of $\Omega$ in the classical sense of quasi-conformal maps, we choose this terminology as they are related (see \cite[Corollary VII.3.4]{gm}).
\end{defn}

Note that by the triangle inequality, for a $K$-quasicircle we have 
\begin{equation}\label{eq:diam_qc}
\abs{w_1-w_2}\le \diam \bb{C_{w_1,w_2}}:= \underset{z,w\in C_{w_1,w_2}}\sup\abs{z-w}\le  \underset{z,w\in C_{w_1,w_2}}\sup\abs{z-w_1}+\abs{w_1-w}< 2K\abs{w_1-w_2}.
\end{equation}

This does not imply that the curve has finite length, e.g., the classical Koch Snowflake is a quasicircle, of infinite length.

The following statement is standard; we include a proof for completeness. 

\begin{claim}\label{claim:cover}
    Let $\gamma$ be a quasicircle with dilation $K$. For every $L>1$ there exists M=M(L, K), depending only on $L$ and $K$ such that for every $\delta\in\bb{\frac{\diam(\gamma)}L,\diam(\gamma)}$, $\gamma$ can be partitioned into $k\leq M$ adjacent  sub-arcs, $\gamma=I_1\cup \dots \cup I_k$, with $$\diam(I_j)=\delta,\ j<k,\quad \diam(I_k)\leq\delta.$$
\end{claim}

\begin{proof}
Let $\gamma: [0,1]\mapsto\mathbb{C}$ be a parameterization of $\gamma$, and define recursively $$t_0=0,\quad t_j=\min\{1, \min\{t>t_{j-1}\,,\,\diam(\gamma([t_{j-1}, t]))=\delta\}\}, \quad k:=\min\{j\,:\,t_j=1\}.$$
Note that the collection $\gamma_j:=\gamma\bb{\sbb{t_{j-1},t_j}}$ satisfies the requirements of the lemma. To conclude the proof it is enough to bound $k$ from above, independently of $\delta$.

Since by definition $\diam(\gamma)\le L\cdot \delta$, one can cover $\gamma$ by at most $\sim K^2L^2$ disks of diameter $\frac\delta {2K}$. 
Assume one of these disks, $B$, contains both $\gamma(t_i)$ and $\gamma(t_j)$  for some $i<j<k$. Then by inclusion
$$\diam(\gamma([t_i,t_j]))\geq\diam(\gamma([t_i, t_{i+1}]))=\delta.$$ However, if $\gamma(t_j),\gamma(t_i)\in B$, then since $\gamma$ is a quasicircle, \eqref{eq:diam_qc} implies
$$\diam(\gamma([t_i,t_j]))<2 K|\gamma(t_i)-\gamma(t_j)|<2K\diam(B)<2K\frac{\delta}{2K}=\frac{\delta}{2}$$
which leads to a contradiction. We conclude that every disk contains at most 2 points in  $\bset{\gamma(t_j)}$ implying that $k\lesssim  K^2L^2$.

\end{proof}

\begin{defn}\label{def:doubling_c}
We say a measure, $\nu$, supported on a Jordan arc, is \emph{doubling} if there exists a constant $C>0$ such that for every adjacent sub-arcs, $I,J$, (i.e., $\overline I\cap\overline J\neq\emptyset$) if $\diam(I)\le \diam(J)$ then $\mu(I)\le C\mu(J)$ (see \cite[p.244]{gm}).
\end{defn}

\begin{defn}\label{def:doubling_m}
We say a measure defined on a metric space, $\mu$, is \emph{doubling} if there exists a constant $C>1$ so that for every disk, $B$, $\mu(2B)\le C\mu(B)$.
\end{defn}

In quasicircles, the two notions of doubling coincide, as the following lemma shows.
\begin{lem}\label{lem:doubling_measure_quasidisk}
Let $\mu$ be a measure defined on a quasicircle, $\partial\Omega$. Then $\mu$ is doubling in the sense  of Definition \ref{def:doubling_c} if and only if doubling in the sense  of Definition \ref{def:doubling_m}.
\end{lem}
\begin{proof}
First, assume that $\mu$ satisfies Definition \ref{def:doubling_c}, and let $B$ be a disk centered at $\partial\Omega$, and $\gamma$ be the smallest connected arc of $\partial\Omega$ containing $2B\cap\partial\Omega$.  The distance between endpoints of $\gamma$ is bounded by $2\diam(2B)=4\diam(B)$, implying that the diameter of $\gamma$ is bounded by $8K_\Omega\diam(B)$.
 
 Let $I$ be the arc of $B\cap\partial\Omega$ that contains the center of $B$. Observe that $\diam(I)\geq\frac12\diam(B)$, since it contains the center of $B$. Using Claim \ref{claim:cover}, we see that each one of the connected components of $\gamma\setminus I$ can be covered by at most $M=M(K,K)$ adjacent arcs $I_1,\dots,I_k$ with $\diam(I_j)=\diam(B)$, $\diam(I_k)\leq\diam(B)$. Therefore, by the doubling property of $\mu$, $\mu(I_j)\leq C^j\mu(I)$, implying that
 $$\mu(2B)\leq\mu(\gamma)\leq \mu(I)\left(1+2\sum_{j=1}^M C^j\right)\leq C_1 \mu(B).$$  

 Now, assume that $\mu$ satisfies Definition \ref{def:doubling_m} and let $I\subset\partial\Omega$ be an arc. Let $z_1, z_2\in I$ be two points satisfying $\abs{z_1-z_2}>\frac{\diam (I)}2$, and let $z_0$ be a point on $C_{z_1, z_2}$ which is equidistant from $z_1$ and $z_2$. We will show that 
 $$B\left(z_0,\frac{\diam(I)}{8K_\Omega}\right)\cap\partial\Omega\subset I.$$

 Indeed, if $z\in\partial\Omega\setminus I$ then since either $z_1\in C_{z,z_0}$ or $z_2\in C_{z,z_0}$ we see that
 $$|z-z_0|\ge\frac{\diam(C_{z,z_0})}{K_\Omega}\ge\frac{\min\{\abs{z_1-z_0},\abs{z_2-z_0}\}}{K_\Omega}\ge\frac{\diam(I)}{8K_\Omega}.$$

Let $J$ be a curve adjacent to $I$ which satisfies $\diam(J)\le\diam(I)$. In particular, $J\subset B\left(z_0,2\diam(I)\right)=B\left(z_0,16K_\Omega\cdot \frac{\diam(I)}{8K_\Omega}\right)$ implying that

 $$\mu(J)\le\mu\bb{B\left(z_0,16K_\Omega\cdot \frac{\diam(I)}{8K_\Omega}\right)}\le C_1\mu\bb{B\left(z_0,\frac{\diam(I)}{8K_\Omega}\right)}\le C_1\mu(I)$$
 where $C_1$ depends only on $K_\Omega$ and $C$ in Definition \ref{def:doubling_m}.
\end{proof}

It is known that if $\Omega$ is a quasidisk, then its corresponding harmonic measure is doubling in the sense of Definition \ref{def:doubling_c} (see \cite[Theorem VII.3.5]{gm}). Thus, the harmonic measure on the quasicircle satisfies Lemma \ref{lem:doubling_measure_quasidisk}.

We will also use the following standard property of quasicircles;

\begin{claim}\label{clm:quasi_finite_curves}
Let $\Omega\subset\C$ be a quasidisk and denote by $K_\Omega$ the dilatation of $\Omega$. For every disk, $B$, for every $L>1$ and for every pair, $C_1,C_2$, of connected components of $\partial\Omega\cap L\cdot B$, which intersect $B$
$$
\frac{(L-1)\diam(B)}{4K_\Omega}\le \dist(C_1\cap B,C_2\cap B)\le \diam(B).
$$
In particular, the number of connected components of $L\cdot B\cap\partial\Omega$ that intersect $B$ is bounded by a constant depending on $K_\Omega$ and $L$ alone. 
\end{claim}
\begin{proof}
Fix a disk, $B$, and let  $C_1,C_2$ be a pair of connected components of $\partial\Omega\cap L\cdot B$ which intersect $B$. By definition of quasidisks, for every $x\in C_1, y\in C_2$ there exists $C_{x,y}\subset\partial\Omega\cap L\cdot B$ satisfying
$$
\diam(C_{x,y})\le 2K_\Omega\abs{x-y}.
$$
Since $C_1,C_2$ are disjoint we can write $C_{x,y}\cap L\cdot B=\bb{C_1\cap C_{x,y}}\uplus  \bb{C_2\cap  C_{x,y}}$ . Each of those arcs, $C_j\cap C_{x,y}$, intersects both $\partial B$ and $\partial\bb{L\cdot B}$. 
We conclude that 
$$
\underset{x\in C_1\cap B\atop y\in C_2\cap B}\inf\;\abs{x-y}\ge \frac1{2K_\Omega}\underset{x\in C_1\cap B\atop y\in C_2\cap B}\inf\;\diam( C_{x,y})\ge \frac{\dist\bb{\partial B,\partial \bb{L\cdot B}}}{2K_\Omega}=\frac{(L-1)\diam(B)}{4K_\Omega}.
$$

The final assertion, regarding the number of connected components, follows from the lower bound. If the distance between any two components within $B$ is at least $(L-1)\frac{\diam(B)}{4K_\Omega}$, then the number of such components must be bounded by 
$$
\frac{\diam(B)^2}{\dist(C_1\cap B,C_2\cap B)^2}\lesssim \frac{K_\Omega^2}{(L-1)^2},
$$
concluding the proof.
\end{proof}

\subsection{Auxiliary Results for Rotation}
In this subsection we reveal intriguing properties of the rotation. The first lemma shows that one can estimate the rotation using integration over curves in $\Omega$:
\begin{lem}\label{lem:rotation_integral}
Fix $w\in\partial\Omega$ and $\delta>0$, and denote by $\Omega_{w,\delta}$ the connected component of $\Omega\setminus B(w,\delta)$ containing $z_0$. Then, every curve, $C_y\subset\Omega_{w,\delta}$, connecting $y\in\partial B(w,\delta)\cap\Omega_{w,\delta}$ with $z_0$ satisfies
$$
\abs{\im{\integrate{C_y}{}{\frac1{\xi-w}}\xi}-\log\bb{\rot(w,\delta)}}\le  3\pi.
$$
\end{lem}
\begin{proof}
Since the branch of the argument is chosen so that $\argu_{[w]}(z_0-w)\in(-\pi,\pi]$ we get that
$$
\im{\integrate{ C_y}{}{\frac1{\xi-w}}\xi}-\pi\le \argu_{[w]}(y-w)\le \im{\integrate{ C_y}{}{\frac1{\xi-w}}\xi}+\pi.
$$
Next, for two points $x, y\in\partial B(w,\delta)\cap\Omega_{w,\delta}$, let $\sigma\subset\partial B(w,\delta)$ be chosen so that the domain bounded by $ C_x+\sigma- C_y$, which is contained in $\Omega_{w,\delta}$, does not contain $w$. Then, since the mapping $\xi\mapsto\frac1{\xi-w}$ is holomorphic in a neighborhood of this domain and the curve is closed,
\begin{eqnarray*}
0&=&\integrate{ C_x+\sigma- C_y}{}{\frac1{\xi-w}}\xi=\integrate{ C_x}{}{\frac1{\xi-w}}\xi+\integrate{\sigma}{}{\frac1{\xi-w}}\xi-\integrate{ C_y}{}{\frac1{\xi-w}}\xi\\
&\Rightarrow& \im{\integrate{ C_x}{}{\frac1{\xi-w}}\xi}-2\pi\le \im{\integrate{ C_y}{}{\frac1{\xi-w}}\xi}\le  \im{\integrate{ C_x}{}{\frac1{\xi-w}}\xi}+2\pi,
\end{eqnarray*}
as the rotation along $\sigma$ is bounded by the rotation of a circle, which is $2\pi$. Overall, we conclude that
\begin{equation*}
    \abs{\im{\integrate{ C_y}{}{\frac1{\xi-w}}\xi}-\argu_{[w]}(x-w)}\le3\pi.
\end{equation*}
In particular the argument above holds for $y^*\in\partial B(w,\delta)\cap\Omega$ which satisfies $\rot(w,\delta)=\exp\bb{\argu_{[w]}(y^*-w)}$.
\end{proof}

The next lemma demonstrates that, unlike harmonic measure, the rotation remains stable; perturbing the center of a disk, $B(w,\delta)$, by less than its diameter does not significantly alter the rotation.


\begin{lem}\label{lem:rotation_err_in}
For every $\xi,w\in\partial\Omega$ if $\abs{\xi-w}<\delta$, then
$$
\abs{\log\bb{\rot(\xi,\delta))}-\log\bb{\rot(w,\delta))}}\le 10\pi.
$$
\end{lem}

\begin{proof}
Let $\Omega_1$ denote the connected component of $\Omega\setminus \bb{B(\xi,\delta)\cup B(w,\delta)}$ containing $z_0$, and let $ C_w, C_\xi\subset\Omega_1$ be two curves connecting $z_0$ with $\partial B(w,\delta)$ and $\partial B(\xi,z)$ respectively. Let $\sigma\subset\partial B(w,\delta)\cup \partial B(\xi,\delta)$ be so that the domain bounded by $\Gamma:= C_w+\sigma- C_w$ does not contain either points, $\xi$ or $w$. Since $B(w,\delta)\cap B(\xi,\delta)\neq\emptyset$, every curve either circles both points or circles none of them. Now in the domain bounded by $\Gamma$ both functions $z\mapsto\frac1{z-\xi}$ and $z\mapsto\frac1{z-w}$ are holomorphic and therefore,
$$
\underset\Gamma\int \frac1{z-\xi}dz=0=\underset\Gamma\int \frac1{z-w}dz,
$$
and in particular
$$
\abs{\im{\underset{ C_\xi}\int \frac1{z-w}dz}-\im {\underset{ C_w}\int \frac1{z-w}dz}}\le\abs{\im{\underset{\sigma}\int \frac1{z-w}dz}}\le 4\pi,
$$
The last piece of the puzzle we need is to observe that since $ C_w, C_\xi\subset\Omega_1$ then the number of times each of these curves circles $\xi$ has to be equal to the number of times it circles $w$ for otherwise, the curve separates between the two points, which is impossible by the way $\Omega_1$ was defined. Note that since the winding index of $C_\xi$ around $\xi$ and the winding index of $C_\xi$ around $w$ is the same,
$$
\frac1{2\pi i}\underset{ C_\xi}\int \frac1{z-w}dz=\frac1{2\pi i}  \underset{ C_\xi}\int \frac1{z-\xi}dz.
$$
Overall, using Lemma \ref{lem:rotation_integral} and the estimates above, we see that
\begin{align*}
&\abs{\log\bb{\rot(B(\xi,\delta))}-\log\bb{\rot(B(w,\delta))}}\le 6\pi+\abs{\im{\integrate{ C_\xi}{}{\frac1{z-\xi}}z}-\im{\integrate{ C_w}{}{\frac1{z-w}}z}}\\
&\quad\quad\quad\le 6\pi+\abs{\im{\integrate{ C_\xi}{}{\frac1{z-\xi}}z}-\im{\integrate{ C_\xi}{}{\frac1{z-w}}z}}+\abs{\im{\integrate{ C_\xi}{}{\frac1{z-w}}z}-\im{\integrate{ C_w}{}{\frac1{z-w}}z}}\le 10\pi.
\end{align*}
\end{proof}

\begin{rmk}\label{rmk:rot_crosscut}
The same proof shows that if two disks, $B_1,\;B_2$ are tangential or intersecting and the connected component of $\Omega\setminus(B_1\cup B_2)$ containing $z_0$ is contained in the connected component of $\Omega\setminus B_j$ containing $z_0$ for $j\in\bset{1,2}$ then their rotations are equivalent up to numerical constants. 

In particular, for every crosscut,  $C$, and every $\xi\in \nu_C$, $\rot(C)\sim \rot(\xi,\dist(\xi,C)+\diam(C))$. To see this let $w_\xi\in C$ satisfy $\abs{w_\xi-\xi}=\dist(\xi,C)$, then the disks $B(\xi,\abs{\xi-w_\xi}+\diam(C))$ and $B(w_\xi,\diam(C))$ are tangential and, since $B(w_\xi,\diam(C))\subset B(\xi,\dist(\xi,C)+\diam(C))$, the corresponding connected components are nested as required. Following the remark above $\rot(C)\sim\rot(w_\xi,\diam(C))\sim \rot(\xi,\dist(\xi,C)+\diam(C))$.
\end{rmk}
Next, it would be interesting to keep the disk centred at $\zeta$ and perturb the diameter. However, in general, this cannot be bounded as the domain may rotate a lot in any annulus. Nevertheless, for quasidisks, this can be uniformly bounded, which is the topic of the next lemma.

\begin{figure}[htp]
\begin{minipage}[c]{.4\textwidth}
\centering
\includegraphics[scale=1.]{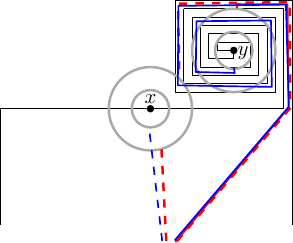}
\label{fig:ro_pert}

\end{minipage}\hfill
\begin{minipage}[c]{.6 \textwidth}
\caption{This is an example showing how to disks of the same radius, centered at $x$ and $y$, could be close but have very different rotations. The same is true for two concentric disks (centered at $x$) with different radii. The dashed curves are curves connecting $z_0$ with disks centered at $x$. The red lines are curves connecting $z_0$ with larger disks, while the blue ones connect it to the smaller disks.}

\end{minipage}
\end{figure}


\begin{lem}\label{lem:rot_quasidisk_radius}
Let $\Omega\subset\C$ be a quasidisk of dilatation $K_\Omega$, fix a point $w\in\partial\Omega$, and let $\delta_2<\delta_1$. Then 
$$
\abs{\log\bb{\rot(B_2)}-\log\bb{\rot(B_1)}}\le 120  K_\Omega\frac{\log\bb{\frac1{\delta_2}}}{\log\bb{\frac1{\delta_1}}} \quad\quad\quad\quad ,\text{ for }B_j=B\bb{w,\delta_j}.
$$
\end{lem}
The same proof can be used to show that John domains also satisfy this property, with some constant $K=K(\Omega)$.
\begin{proof}
Let $C\subset\Omega\setminus B_2$ be a curve connecting $z_0$ and $p\in\partial B_2$, and let $C'\subset\Omega\setminus B_1$ is the sub-arc of $C$ connecting $z_0$ to $q\in\partial B_1$. Following Lemma \ref{lem:rotation_integral}
$$
\abs{\im{\integrate{C}{}{\frac1{\xi-w}}\xi}-\log\bb{\rot(B_2)}}\le  3\pi;\quad \abs{\im{\integrate{C'}{}{\frac1{\xi-w}}\xi}-\log\bb{\rot(B_1)}}\le  3\pi,
$$ 
implying that
\begin{equation*}
\abs{\log\bb{\rot(B_2)}-\log\bb{\rot(B_1)}}\le 6\pi +\abs{\im{\integrate{C\setminus C'}{}{\frac1{\xi-w}}\xi}}.
\end{equation*}

Define $L=[p,q]$. If $L\subset\Omega$, we define $\Gamma=C\setminus C'+L$ and note that 
$$
\abs{\im{\integrate{C\setminus C'}{}{\frac1{\xi-w}}\xi}}=\abs{\im{\integrate{L}{}{\frac1{\xi-w}}\xi}}\le\pi,
$$
since the rotation of a straight line is just a change in angle, which is bounded by $\pi$.

If $L\not\subset\Omega$ then since the initial point of $L$ and the terminal point of $L$ both lie in $\Omega$, it must intersect the boundary of $\Omega$ an even number of times. Denote the points where it leaves $\Omega$ by $\bset {b_\nu}$ and the ones it re-enters $\Omega$ by $\bset{e_\nu}$, i.e., the segment $\bb{e_\nu,b_{\nu+1}}\subset\Omega$ while $\sbb{b_\nu,e_\nu}\subset\C\setminus\Omega$. Since $\Omega$ is a quasidisk, for every $\nu$
$$
\ell_\Omega(b_\nu,e_\nu)\le K_\Omega \ell_{\C\setminus\Omega}(b_\nu,e_\nu)=K_\Omega\cdot \abs{e_\nu-b_\nu},
$$
where $\ell_\Omega$ denotes the distance inside $\Omega$. 
\vspace{-0.3em}
\begin{figure}[htp]
\begin{minipage}[c]{.505\textwidth}
For every $\nu$, we replace $[b_\nu,e_\nu]\subset L\cap\C\setminus\Omega$ with a curve $C_\nu\subset\Omega$ of diameter bounded by $2K_\Omega\cdot \abs{e_\nu-b_\nu}$. If $C_\nu\cap B_2\neq\emptyset$, we replace the rest of the line, $L$, with the connected component of $C_\nu\setminus B_2$ containing $b_\nu$ and stop. We may therefore assume that $C_\nu\cap B_2=\emptyset$. See Figure \ref{fig:c_nu}.

Next, note that the winding number of $C_\nu+\sbb{b_\nu,e_\nu}$ in absolute value is at most $2\pi$ for each point outside domain bounded by this curve. This implies

\end{minipage}\hfill
\begin{minipage}[c]{.47\textwidth}
\centering
\includegraphics[scale=0.85]{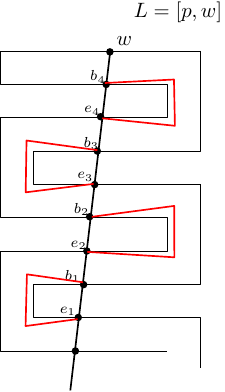}
\caption{The red curves are the curves $C_\nu$ which replace the parts of $C+L$ which is outside $\Omega$.}
\label{fig:c_nu}
\end{minipage}
\end{figure}

$$
\abs{\im{\integrate{C_\nu}{}{\frac1{\xi-w}}\xi}}\le \abs{\im{\integrate{C_\nu+\sbb{b_\nu,e_\nu}}{}{\frac1{\xi-w}}\xi}}+ \abs{\im{\integrate{\sbb{b_\nu,e_\nu}}{}{\frac1{\xi-w}}\xi}}\le3\pi,
$$
as this integral over a line segment is bounded by $\pi$.

Combining these, we see that
\begin{align*}
&\abs{\log\bb{\rot(B_2)}-\log\bb{\rot(B_1)}}\le  6\pi +\abs{\im{\integrate{C\setminus C'}{}{\frac1{\xi-w}}\xi}}\\
&\quad\quad\quad\quad\le 6\pi+\sumit \nu 1N\bb{\abs{\im{\integrate{[e_\nu,b_\nu]}{}{\frac1{\xi-w}}\xi}}+\abs{\im{\integrate{C_\nu}{}{\frac1{\xi-w}}\xi}}}=6\pi+3\pi\#\bset{\nu, C_\nu \text{ coils around }B_2}.
\end{align*}
To conclude the proof, we need to bound the number of curves that coil around $B_2$, which is where the fact that $\Omega$ is a quasidisk is used.

Recall that $C_\nu$ was chosen so that $\diam(C_\nu)\le 2K_\Omega\cdot \abs{e_\nu-b_\nu}$. Then, if $C_\nu$ coils around $B_2$ we obtain
$$
2K_\Omega\abs{e_\nu-b_\nu}\ge \diam(C_\nu)\ge 2 \abs{e_\nu-w}\Rightarrow \abs{e_\nu-b_\nu}\ge\frac{\abs{e_\nu-w}}{K_\Omega}.
$$
Assume that the sequence $\bset{e_\nu}$ is ordered in decreasing distance from $w$, i.e., $\abs{e_\nu-w}>\abs{e_{\nu+1}-w}$ for all $\nu$, let $\nu_k$ denote the sequence defined by
$$
\nu_k:=\min\bset{\nu, \abs{e_\nu-w}<2^{-k+1}\delta_1}.
$$
Note that $\nu_1=1$ as all $\nu$ satisfy $\abs{e_\nu-w}<\delta_1$. While $\bset{\nu_k}$ is monotone increasing, it could be the case that $\nu_k=\nu_{k+1}$. In addition, for every $\nu_k\le \nu<\nu_{k+1}$ (when this interval is non-trivial) we have $\abs{e_\nu-w}\in\left[2^{-k}\cdot\delta_1,2^{-k+1}\delta_1\right)$.

For every $k$, we define
$$
\tau_k:=\lambda_1\bb{\underset{\nu_k\le\nu<\nu_{k+1}}\bigcup \sbb{b_\nu,e_\nu}}\le 2^{-k+2}\delta_1,
$$
(by inclusion of the intervals $[p,w]\cap B\bb{w,2^{-k+2}\delta_1}$ and disjointness of the curves, $C_\nu$) and 
$$
m_k:=\underset{\nu_k\le \nu<\nu_{k+1}, \atop C_\nu \text{ coils around }B_2}\inf\; \abs{e_\nu-b_\nu}\ge\underset{\nu_k\le \nu<\nu_{k+1},\atop C_\nu \text{ coils around }B_2}\inf\;\frac{\abs{e_\nu-w}}{K_\Omega}\ge \frac{\delta_1}{2^k K_\Omega}.
$$
Then, the number of indices $\nu_k\le\nu<\nu_{k+1}$ coiling around $B_2$ is bounded by the total length of $\underset{\nu_k\le\nu<\nu_{k+1}}\bigcup \sbb{b_\nu,e_\nu}$ divided by their minimal length, which is bounded from below by
$$
\frac{\tau_k}{m_k}\le \frac{2^{-k+2}\delta_1}{\frac{\delta_1}{2^k K_\Omega}}=4K_\Omega.
$$
Lastly, the largest index $K$ to consider is when $2^{-K+1}\delta_1\le\delta_2$, i.e.,
$$
K:=\left\lceil1+\frac{\log\bb{\frac1{\delta_2}}}{\log\bb{\frac1{\delta_1}}\cdot\log(2)}\right\rceil\le 4\frac{\log\bb{\frac1{\delta_2}}}{\log\bb{\frac1{\delta_1}}}.
$$
Overall, we conclude that
\begin{align*}
\abs{\log\bb{\rot(B_2)}-\log\bb{\rot(B_1)}}&\le  \cdots\le 6\pi+2\pi\#\bset{\nu, C_\nu\text{ coils around }B_2}\le \\
&\le  6\pi+8\pi\cdot K_\Omega\cdot K\le 6\pi+32\pi\cdot K_\Omega\frac{\log\bb{\frac1{\delta_2}}}{\log\bb{\frac1{\delta_1}}}\le 120  K_\Omega\frac{\log\bb{\frac1{\delta_2}}}{\log\bb{\frac1{\delta_1}}},
\end{align*}
concluding the proof.
\end{proof}


The next lemma we present in this auxiliary subsection is kind of a mean-value theorem for the argument of the derivative of holomorphic functions. To prove the lemma, we will require an observation about the diameter of images of arcs.
\begin{obs}\label{obs:distortion_vs_curve}
Let $\phi:\D\rightarrow\C$ be a conformal map and fix $I\subset\partial \D$ some arc with $\lambda_1(I)<\frac12$. Let  $z_I$ be the point representing $I$, and let $C_I$ be a crosscut of $I$. Then
\begin{enumerate}
\item $\dist(\phi(z_I),\nu_I)\dless{\Omega} \diam(C_I)$.
\item If $\Omega$ is a quasidisk, then $\diam(C_I)\dsim{K_\Omega}\bb{1-\abs{z_I}}\abs{\phi'(z_I)}$, where $K_\Omega$ is the dilatation of $\Omega$.
\end{enumerate}
\end{obs}
See \cite[Ex 8, p.153]{gm} for the first part, and on \cite[p. 261]{gm} for the second part. 

\begin{lem}\label{lem:rotation_differentiation} Let $\phi:\D\rightarrow\Omega$ be a conformal map. For every $I\subset\T$ satisfying $\lambda_1(I)\le \frac12$ there exists $\eta\in I$ for which
$$
\abs{\argu_{[\phi(\eta)]}\sbb{\frac{\phi(z_I)-\phi(\eta)}{z_I-\eta}}-\argu\sbb{\phi'(z_I)}}\le c_1+c_2(\abs{z_I})\log\bb{\frac1{1-\abs{z_I}}},
$$
where
\begin{itemize}
    \item $z_I$ is the point representing $I$.
    \item $\argu\phi'(z)$ is normalized so that $\argu\phi'(0)\in[-\pi,\pi)$.
    \item $\argu_{[\phi(\eta)]}$ is normalized so that the function $(z,w)\mapsto \argu_{[\phi(\eta)]}\bb{\frac{\phi(z)-\phi(w)}{z-w}}$ agrees with the branch of the argument  satisfying  $\argu_{[\phi(\eta)]}\bb{\frac{\phi(0)-\phi(\eta)}{\eta-0}}\in(-\pi,\pi]$ (in particular, when $z=w$, $\argu_{[\phi(\eta)]}(\phi'(z))$ is well-defined).
    \item $c_1=c_1(\Omega)>0$ is a constant, and $c_2(r)=c_2\bb{r,\diam\bb{\Omega}}\searrow 0$.
\end{itemize}
If $\Omega$ is an $\alpha$-H\"older domain, then $c_1$ depends on $\alpha$ alone, and $c_2\equiv 0$.
\end{lem}
\begin{rmk}
Every quasicircle is $\alpha$-H\"older for some $\alpha$ depending on $K_\Omega$ (\cite{pombook}, Theorem 5.2). In addition, in this case the constant $c_1$ depends on the dilatation alone instead of on the set. We conclude that statement holds for any quasidisk with some $c_1=C(K_\Omega)$ and $c_2\equiv 0$.  
\end{rmk}
The proof relies heavily on ideas from the proof of McMillan's twist theorem (see, e.g., \cite[p.142]{pombook}).

\begin{proof}
We will first show that there exists $\eta\in I$ for which
\begin{equation}\label{eq:diff_rot_1}
\abs{\argu_{[\phi(\eta)]}\sbb{\frac{\phi(z_I)-\phi(\eta)}{z_I-\eta}}-\argu_{[\phi(\eta)]}\sbb{\phi'(z_I)}}\dless \Omega 1.
\end{equation}

Following \cite[Lemma 6.19]{pombook}, with the point $z=z_I$ and the interval $I$, there exists a point $\eta\in I$ such that the non-euclidean segment $A$, connecting $z_I$ and $\eta$ , satisfies
\begin{equation}\label{eq:close_pt}
\abs{\phi(z_I)-\phi(\eta)}\le \lambda_1(\phi(A))\lesssim \dist\bb{\phi(z_I),\partial\Omega}\lesssim \diam(\phi(C_I)),
\end{equation}
where $C_I$ is a crosscut of $I$, following Observation \ref{obs:distortion_vs_curve}. 

Using the point $\eta$, we define the map $\psi:\D\times\D\rightarrow\R$ by
$$
\psi(z,w):=\begin{cases}
			\argu_{[\phi(\eta)]}\bb{\frac{\phi(z)-\phi(w)}{z-w}}&, z\neq w\\
			\argu_{[\phi(\eta)]}\bb{\phi'(z)}&, z= w
			\end{cases}.
$$
This map is continuous on $\D\times\D$. In addition, by the triangle inequality
$$
\abs{\psi(z_I,z_I)-\psi(\eta,z_I)}\le\abs{\psi(z_I,z_I)-\psi(\zeta,z_I)}+\abs{\psi(\zeta,z_I)-\psi(\eta,z_I)}=S_1+S_2,
$$
where $\eta$ is the end-point of $A$, and $\zeta\in A$ is a point satisfying $\rho(\zeta,z_I)=1$. 

To bound $S_1$, we use \cite[Ex 4, p.13]{pombook},
$$
S_1=\abs{\psi(z_I,z_I)-\psi(\zeta,z_I)}=\abs{\argu_{[\phi(\eta)]}\bb{\phi'(z_I)}-\argu_{[\phi(\eta)]}\bb{\frac{\phi(\zeta)-\phi(z_I)}{\zeta-z_I}}}\le8\rho(z_I,\zeta)+\frac\pi2\le 10.
$$

To bound $S_2$ we will need to work harder. Define the non-euclidean sub-segment $\tilde A:=\bset{w\in A, \rho(w,z_I)\ge1}$. Then 
\begin{align*}
S_2=\abs{\psi(\zeta,z_I)-\psi(\eta,z_I)}=\underset{\tilde A}\int\;d\abs{\argu_{[\phi(\eta)]}\sbb{\phi(z)-\phi(z_I)}}\le\integrate{\tilde A}{}{\frac{\abs{\phi'(z)}}{\abs{\phi(z)-\phi(z_I)}}}{\abs z}\lesssim\frac{1}{\dist(\phi(z_I),\partial\Omega)} \integrate{\tilde A}{}{\abs{\phi'(z)}}{\abs z}\nonumber\\
\le\frac{1}{\dist(\phi(z_I),\partial\Omega)} \integrate{A}{}{\abs{\phi'(z)}}{\abs z}\lesssim \frac{1}{\dist(\phi(z_I),\partial\Omega)}\cdot \dist(\phi(z_I),\partial\Omega)\sim 1,
\end{align*}
using the fact that following \cite[Cor. 1.5]{pombook}, and Koebe's distortion theorem, for every $z\in\tilde A$
$$
\abs{\phi(z)-\phi(z_I)}\ge \abs{\phi'(z_I)}\bb{1-\abs{z_I}^2}\frac{\tanh(\rho(z,z_I))}4\ge \dist(\phi(z_I),\partial\Omega)\cdot \frac{\tanh(1)}4\sim\dist(\phi(z_I),\partial\Omega).
$$
Combining both estimates together, we see that
$$
\abs{\argu_{[\phi(\eta)]}\sbb{\frac{\phi(z_I)-\phi(\eta)}{z_I-\eta}}-\argu_{[\phi(\eta)]}\sbb{\phi'(z_I)}}=\abs{\psi(\eta,z_I)-\psi(z_I,z_I)}\le S_1+S_2\lesssim 1,
$$
where the constant depends on $\Omega$, and if $\Omega$ is a quasidisk, then they depend on the dilatation, $K_\Omega$, alone.
This concludes the proof of \eqref{eq:diff_rot_1}. 

To conclude the proof of the lemma, we need to compare $\argu_{[\phi(\eta)]}\sbb{\phi'(z_I)}$ and $\argu\sbb{\phi'(z_I)}$. 

We will use the branch of $\argu\bb{\frac{\phi(z)-\phi(w)}{z-w}}$ that agrees with $\argu(\phi'(z))$ when $z=w$, i.e.,  $\argu\phi'(0)\in(-\pi,\pi]$. The function $(z,w)\mapsto \argu_{[\phi(\eta)]}\bb{\frac{\phi(z)-\phi(w)}{z-w}}-\argu\bb{\frac{\phi(z)-\phi(w)}{z-w}}$ is a constant function as a difference between two branches of the argument. Then, evaluating the difference at $\phi'(z_I)$ is the same as evaluating it at $\frac{\phi(0)-\phi(\eta)}{0-\eta}$. Namely,
\begin{align*}
\abs{\argu_{[\phi(\eta)]}\sbb{\phi'(z_I)}-\argu\sbb{\phi'(z_I)}}&=\abs{\argu_{[\phi(\eta)]}\bb{\frac{\phi(0)-\phi(\eta)}{0-\eta}}-\argu\bb{\frac{\phi(0)-\phi(\eta)}{0-\eta}}}\\
&\le \pi+\abs{\argu\bb{\frac{\phi(0)-\phi(\eta)}{0-\eta}}}\le \pi+\abs{\argu\bb{\frac{\phi(z_I)-\phi(0)}{z_I}\cdot \frac{z_I}\eta\cdot \frac{\phi(\eta)-\phi(0)}{\phi(z_I)-\phi(0)}}}\\
&\le 7\pi + \abs{\argu\bb{\frac{\phi(z_I)-\phi(0)}{z_I}}}+\abs{\argu\bb{\frac{\phi(\eta)-\phi(0)}{\phi(z_I)-\phi(0)}}}\le8\pi +\abs{\argu\bb{\frac{\phi(z_I)-\phi(0)}{z_I}}},
\end{align*}
since we assumed that $\argu_{[\phi(\eta)]}\bb{\frac{\phi(0)-\phi(\eta)}{\eta-0}}\in(-\pi,\pi]$ and  
$$
\abs{\argu\bb{\frac{\phi(\eta)-\phi(0)}{\phi(z_I)-\phi(0)}}}=\abs{\argu\bb{1+\frac{\phi(z_I)-\phi(\eta)}{\phi(0)-\phi(z_I)}}}\le\pi
$$
as $\abs{\phi(\eta)-\phi(z_I)}\ll\abs{\phi(0)-\phi(z_I)}$ following \eqref{eq:close_pt}, therefore $\argu\bb{\frac{\phi(\eta)-\phi(0)}{\phi(z_I)-\phi(0)}}$ is bounded and well defined.

Lastly, to bound $\abs{\argu\bb{\frac{\phi(z_I)-\phi(0)}{z_I}}}$, recall that for every brunch of the argument, there exists a brunch of the logarithm satisfying $\argu[z]=\Im\log(z)$. Then
\begin{align*}
\abs{\argu\bb{\frac{\phi(z_I)-\phi(0)}{z_I}}}&=c_1+\abs{\int_{\sbb{\frac{z_I}2,z_I}} d\arg\bb{\frac{\phi(w)}w}}\le  c_1+\int_{\sbb{\frac{z_I}2,z_I}}\abs{d\arg\bb{\frac{\phi(w)}w}}\\
&\le c_1+\int_{\sbb{\frac{z_I}2,z_I}}\abs{\frac d{dw}\bb{\log\bb{\phi(w)}-\log(w)}}=c_1+\int_{\sbb{\frac{z_I}2,z_I}}\frac{\abs{\phi'(w)}}{\abs{\phi(w)}}\abs{dw}+\int_{\sbb{\frac{z_I}2,z_I}}\frac{1}{\abs{w}}\abs{dw}.
\end{align*}

Now, for any simply connected bounded domain, $\Omega$, by Koebe's distortion theorem,
$$
\int_{\sbb{\frac{z_I}2,z_I}}\frac{\abs{\phi'(w)}}{\abs{\phi(w)}}\abs{dw}\le \int_{\sbb{\frac{z_I}2,z_I}}\frac{\dist(\phi(w),\partial\Omega)}{1-\abs{w}}\abs{dw} \leq c_2(\abs{z_I})\log\bb{\frac1{1-\abs{z_I}}}.
$$
If $\Omega$ is a quasidisk, then following Observation \ref{obs:distortion_vs_curve} part 2
$$
\int_{\sbb{\frac{z_I}2,z_I}}\frac{\abs{\phi'(w)}}{\abs{\phi(w)}}\abs{dw}\le \int_{\sbb{\frac{z_I}2,z_I}}\frac{\dist(\phi(w),\partial\Omega)}{1-\abs{w}}\abs{dw} \overset{K_\Omega}\lesssim\int_{\sbb{\frac{z_I}2,z_I}}\frac{\bb{1-\abs{w}}\abs{\phi'(w)}}{1-\abs{w}}\abs{dw} \le \diam(\Omega).
$$
If $\Omega$ is $\alpha$-H\"older, then
$$
\int_{\sbb{\frac{z_I}2,z_I}}\frac{\abs{\phi'(w)}}{\abs{\phi(w)}}\abs{dw}\le \int_{\sbb{\frac{z_I}2,z_I}}(1-\abs w )^{-\alpha}\abs{dz}\leq c_1,
$$
and $c_1$ depends only on $\alpha$, concluding the proof of the lemma.
\end{proof}

%% file: sections/relations.tex
In this section, we prove three key 'folklore results': Theorem \ref{thm:fine=distortion_quasidisks}, Theorem \ref{thm:relations}, and the Main Lemma, Lemma \ref{lem:main}. We bundled all three proofs together not only as they rely on one another, they also share similar ideas and intricacies.
\subsection{Proof of The Main Lemma.}
For the reader's convenience, we restate the lemma \ref{lem:main} here:
\begin{lem*}[The Main Lemma]
Let $\Omega\subset\C$ be a simply connected domain and let $\phi:\D\rightarrow\Omega$ be a Riemann map sending $0$ to $z_0$. Then
\begin{itemize}
\item[]
\vspace{-1cm}
\begin{equation*}
\hspace{-6cm}\tag{HM*}
\abs{\phi'\bb{z_I}} \sim\frac{\dist(\phi(z_I),\partial\Omega)}{\lambda_1(I)}\sim \frac{\diam(C_I)}{\lambda_1(I)}\lesssim \frac{\diam(\nu_I)}{\lambda_1(I)}.
\end{equation*}
\item []
\vspace{-1cm}
\begin{equation*}
\hspace{-10.5cm}\tag{R*}  \lim_{r\to1^-}\frac{\log\abs{\phi'^{-i}\bb{z_I}}}{\log\rot(C_I)}=1.
\end{equation*}
\item [] 
\vspace{-1cm}
 \begin{equation*}
\hspace{-9.5cm} \tag{D*} \rho(z_B,z_{B'})\sim \frac{\omega\bb{C_{I_{B,B'}}}}{\min\bset{\omega(B),\omega(B')}},
\end{equation*}
for every disjoint disks $B,B'$ centered at $\partial\Omega$, and $I_{B,B'}=\sbb{\frac{z_B}{\abs{z_B}},\frac{z_{B'}}{\abs{z_{B'}}}}\subset\partial\D.$
\end{itemize}
If $\Omega$ is a quasidisk, those properties become stronger. Namely, if $K_\Omega$ is the dilatation of $\Omega$, then
\begin{itemize}[]
\addtolength{\itemindent}{1.2cm}
\item[(HM*)]  $\abs{\phi'\bb{z_I}}\dsim{K_\Omega} \frac{\diam(\nu_I)}{\lambda_1(I)}$.
\item [(R*) ~ ]  $ \abs{\phi'^{-i}\bb{z_I}}\dsim{K_\Omega} \rot(C_I)$.
\item [(D*) ~ ]  If $2B\cap 2B'=\emptyset$, then $\rho(z_B,z_{B'})\overset{K_\Omega}\gtrsim1$.
\end{itemize}
\end{lem*}
\begin{proof}
\eqref{eq:arc_measure} is a corollary of Koebe's distortion theorem:
$$
\dist(\phi(z_I),\partial\Omega)\sim \abs{\phi'(z_I)}\cdot \dist(z_I,\partial\D)=\abs{\phi'(z_I)}\cdot \bb{1-\abs{z_I}}=\abs{\phi'(z_I)}\cdot \lambda_1(I)
$$
where the constant in $\sim$ above is 4 (and $\frac14$ respectively). To conclude the proof, we only need to use Observation \ref{obs:distortion_vs_curve} to see that
$$\dist(\phi(z_I),\partial\Omega)\le C_{abs} \diam(C_I)\le 2C_{abs}\diam(\nu_I),$$
where $C_{abs}$ is some numerical constant concluding the proof of the general case. For quasidisks, the opposite inequality holds by the second half of Observation \ref{obs:distortion_vs_curve} (with a constant which depends on $K_\Omega$) combined with the definition of a quasicircle.

To prove \eqref{eq:arc_rotation}, we turn to Lemma \ref{lem:rotation_differentiation}. Let $\eta\in I$ be a point satisfying
$$
\abs{\argu_{[\phi(\eta)]}\sbb{\frac{\phi(z_I)-\phi(\eta)}{z_I-\eta}}-\argu\sbb{\phi'(z_I)}}\le C_\Omega.
$$
We need to relate the argument $\argu_{[\phi(\eta)]}\sbb{\phi(z_I)-\phi(\eta)}$ with $\rot(C_I)$, since it is not necessarily the case that $\phi(z_I)\in\partial B(\phi(\eta),\delta)$.

Following Observation \ref{obs:distortion_vs_curve},
$$
\abs{\phi(z_I)-\phi(\eta)}\le \dist(\phi(z_I),\nu_I)+\diam(\nu_I)\le M\cdot \delta.
$$
This holds for some absolute constant $M$, which if $\Omega$ is a quasidisk, depends on the dilatation alone. 
Then there exists a sequence of intersecting disks $\bset{B_\ell}_{\ell=1}^{2M}$ so that $B_\ell=B(\xi_\ell,\delta)$ with $\xi_1=\phi(\eta)$ while $\phi(z_I)\in\partial B_{2M}$. Let $\sigma_\ell\subset \partial B_\ell\cup\partial B_{\ell+1}$ be so that $\sumit \ell 1 {2M}\sigma_\ell$ is a curve in $\Omega$ connecting $\partial B_1\cap\Omega$ with $\phi(z_I)$, let $ C_{\phi(\eta)}, C_{\phi(z_I)}\subset\Omega$ be two curves connecting $z_0$ with $\partial B_1$ and $\partial B_{2M}$ respectively. Note that the domain bounded by the curves $ C_{\phi(z_I)}- C_{\phi(\eta)}+\sumit \ell 1 {2M}\sigma_\ell$ does not contain the point $\phi(\eta)$. Then, a similar argument to the one presented in Lemma \ref{lem:rotation_err_in} shows that
\begin{align*}
\abs{\log\bb{\rot(\phi(\eta),\delta)}-\argu_{[\phi(\eta)]}\sbb{\phi(z_I)-\phi(\eta)}}&\le\abs{\im{\integrate{ C_{\phi(\eta)}}{}{\frac1{\xi-\phi(\eta)}}\xi}-\im{\integrate{ C_{\phi(\eta)}+\sumit \ell 1 {2M}\sigma_\ell}{}{\frac1{\xi-\phi(\eta)}}\xi}}+6\pi\\
&\le6\pi+\sumit \ell 1 {2M}\abs{\im{\integrate{\sigma_\ell}{}{\frac1{\xi-\phi(\eta)}}\xi}}\le6\pi\bb{M+1},
\end{align*}
since the change in the argument along each $\sigma_\ell$ is bounded by $2\pi$. Next,
\begin{eqnarray*}
\abs{\log\bb{\rot(C_I)}-\argu\sbb{\phi'(z_I)}}&\le&\abs{\log\bb{\rot(C_I)}-\log\bb{\rot(\phi(\eta),\delta)}}\\
&+&\abs{\log\bb{\rot(\phi(\eta),\delta)}-\argu_{[\phi(\eta)]}\sbb{\phi(z_I)-\phi(\eta)}}\\
&+&\abs{\argu_{[\phi(\eta)]}\sbb{\phi(z_I)-\phi(\eta)}-\argu\sbb{\phi'(z_I)}}\dless\Omega 1,
\end{eqnarray*}
as the first summand is bounded by the way rotation of a curve was defined, the second summand is bounded by the computation above, and the third summand is bounded by Lemma \ref{lem:rotation_differentiation}. If $\Omega$ is a quasidisk then the constant depends on the dilatation alone. Overall,
\begin{align*}
\abs{\phi'^{-i}\bb{z_I}}=\exp\bb{Re\bb{-i\log\bb{\phi'(z_I)}}}&=\exp\bb{Re\bb{-i\bb{\log\abs{\phi'(z_I)}+i\argu\sbb{\phi'(z_I)}}}}\\
&=\exp\bb{\argu\sbb{\phi'(z_I)}}\dsim\Omega \rot(C_I),
\end{align*}
The last equivalence is by the last assertion of Lemma \ref{lem:rotation_differentiation}.

Lastly, to prove \eqref{eq:disjoint_dsks}, note that 
$$
\rho\bb{z_B,z_{B'}}\sim\frac{\abs{z_B-z_{B'}}}{\min\bset{1-\abs{z_B},1-\abs{z_{B'}}}}\sim\frac{\omega\bb{C_{I_{B,B'}}}}{\min\bset{\omega(B),\omega(B')}}.
$$
If $\Omega$ is a quasidisk, and $2B\cap 2B'=\emptyset$, then, since the harmonic measure of $\Omega$ is doubling, 
$$
\frac{\omega\bb{C_{I_{B,B'}}}}{\min\bset{\omega(B),\omega(B')}} \dmore{K_\Omega}1
$$
concluding the proof.
\end{proof}

\subsection{Proof of Theorem \ref{thm:relations}}
In light of the main lemma, the proof of this theorem is reduced to making the correct finite to one correspondence between points in $\D$ with certain derivative, and disks centred in $\partial\Omega$ with certain harmonic measure and rotation.

\begin{proof}
Recall that we defined the arc $A_z=\bset{\xi\in\partial\D, \abs{\xi-\frac z{\abs z}}<1-\abs z}$, where $z\in\partial \Omega$. Given a point $z\in\partial\Omega$ we define the corresponding disk $B_z$ to be the smallest disk centered at $\partial\Omega$ containing $\phi(A_z)$. Since $\phi(A_z)\subset\partial\Omega$ and $B_z$ is centered in $\partial\Omega$, then diameter of the disk is between $\diam(\phi(A_z))$ and $2\diam(\phi(A_z))$. By definition, 
$$
\omega(B_z)\ge\omega(\phi(A_z))=\lambda_1(A_z)\sim 1-\abs{z}\dsim{K_\Omega}\frac{\diam(\phi(A_z))}{\abs{\phi'(z)}}
$$
by Observation \ref{obs:distortion_vs_curve}, as $\Omega$ is a quasidisk. However, following Lemma \ref{lem:doubling_measure_quasidisk}, $\omega(B_z)\dsim{K_\Omega}\omega(\phi(A_z))$ as well. In addition,
$$
\rot(B_z)\dsim{K_\Omega} \rot(C_{A_z})\dsim{K_\Omega} \abs{\phi'^{-i}\bb{z}}
$$
following the main lemma, Lemma \ref{lem:main}, and Lemma \ref{lem:rotation_err_in}.

Lastly, let $z,z'\in\D$ be so that $B_z\cap B_{z'}\neq\emptyset$, then following \cite[Corollary 1.5]{pombook}
$$
4\delta=\diam(B)+\diam(B')\ge\abs{\phi(z)-\phi(z')}\ge\frac14\bb{1-\abs z^2}\abs{\phi'(z)}\cdot \rho(z,z')\dsim {K_\Omega}\delta \cdot \alpha(K_\Omega),
$$
which is a contradiction if $\alpha(K_\Omega)$ is chosen large enough.

On the other hand, given a disk, $B$, let $\gamma_B$ be an arc representing $B$, and let $z_B$ be the point representing $B$. Note that since the harmonic measure is doubling (see Lemma \ref{lem:doubling_measure_quasidisk}), by looking at the extension of $\gamma_B$ in $3B$ we may assume that $\diam(\beta_B)\in\bb{\delta,6\delta}$ and
$$
\frac{\omega(B)}{\log^2\bb{\frac1{\omega(B)}}}\dless{K_\Omega} \omega(\beta_B)\dless{K_\Omega} \omega(B).
$$
Define $z_\beta=\xi_\beta\bb{1-\omega(\beta_B)}$, then following the main lemma, Lemma \ref{lem:main}, $z_\beta$ satisfies \eqref{eq:arc_measure} and \eqref{eq:arc_rotation}. Next, note that
$$
\rho(z_B,z_\beta)\sim\frac{\abs{\omega(B)-\omega(\beta)}}{\min\bset{\omega(B),\omega(\beta)}}	\begin{cases}
																						\gtrsim 1\\
																						\lesssim\log^2\bb{\frac1{\omega(B)}}
																						\end{cases}.
$$
Since $\log\phi'$ is a Bloch function, we see that
$$
\abs{\log\bb{\phi'(z_B)}-\log\bb{\phi'(z_\beta)}}\lesssim \rho\bb{z_B,z_\beta}\lesssim \log^2\bb{\frac1{\omega(B)}}.
$$
To conclude the proof, it is enough to note that the rotation of the disk and the rotation of the crosscut of $\phi\inv(\beta)$ are proportional but this holds by Lemma \ref{lem:rotation_err_in}.

Lastly, if $2B\cap 2B'=\emptyset$, them following the Main Lemma, Lemma \ref{lem:main}, part \eqref{eq:disjoint_dsks}, $\rho(z_B,z_{B'})\dmore{K_\Omega}1$ concluding the proof.
%
%

\end{proof}
\vspace{-1cm}\subsection{Proof of Theorem \ref{thm:fine=distortion_quasidisks}}

\paragraph*{Proof that $d_\Omega\ge f_\Omega$:} It is enough to show that for every $\eta>0$ there exists a sequence $\bset{r_k}\nearrow1$ so that for every $k$ large enough,
$$
\frac{\log\bb{\lambda_1\bb{L_{a,b,\eta}(r_k)}}}{\log\bb{\frac1{1-r_k}}}\ge \bb{1-a}f_\Omega\bb{\frac1{1-a},\frac{-b}{1-a}}-1-\eta.
$$

Let $\eta'\in\bb{0,\frac{\eta\cdot \alpha}3}$ and $\bset{\delta_k}$ be so that $\limit k\infty\frac{\log\bb{N(\delta_k,\alpha,\gamma,\eta')}}{\log\bb{\frac1{\delta_k}}}\ge f_\Omega(\alpha,\gamma)-\eta\cdot\alpha$. For every $k$ there exists a collection of pairwise disjoint disks, $\bset{B_j^k}_{j=1}^{N(\delta_k,\alpha,\gamma,\eta')}$, centered at $\partial\Omega$ of radius $\delta_k$ satisfying \eqref{eq:N_measure} and \eqref{eq:N_rotation}. By excluding at most a linear portion of the disks in the collection, we may assume without loss of generality that $3B_j^k\cap 3B_\nu^k=\emptyset$ for every $j\neq\nu$. Following Theorem \ref{thm:relations}, for every $j$ the point $z_j^k:=z_{B_j^k}\in\bb{1-\delta_k^\alpha}\T$ satisfies \eqref{eq:derivative_measure},  \eqref{eq:derivative_rotation}, and \eqref{eq:disjoint_dsks}. 

Note that in this case, if $\delta_k$ is small enough (depending on $\eta'$)
$$
\begin{cases}
	\delta_k^{1-\alpha+2\eta'}\le \frac{\delta_k}{C_\Omega\delta_k^\alpha\log^2\bb{\frac1{\delta_k^\alpha}}}\le\abs{\phi'\bb{z_j^k}} \le C_\Omega \frac{\delta_k}{\delta_k^\alpha}\le \delta_k^{1-\alpha-2\eta'}\\
	\delta_k^{\gamma+2\eta'}\le \frac{\rot\bb{B_j^k}}{C_\Omega}\le \abs{\phi'^{-i}\bb{z}}\le C_\Omega \rot\bb{B_j^k}\le\delta_k^{\gamma-2\eta'}
\end{cases}.
$$
We divide $\bb{1-\delta_k^\alpha}\T$ into $N:=\lceil{\delta_k^{-\alpha}}\rceil$ arcs of equal length, and denote this collection $\mathcal P_k$. Note that for every $j\neq\nu$ we have $\rho(z_j^k,z_\nu^k)\sim 1$. By excluding at most a linear portion of the disks, the points $\bset{z_j^k}$ belong to different arcs in the collection $\mathcal P_k$.

For every $r$ we set $r_k=1-\delta_k^\alpha$, for $\delta_k$ small enough so that $r_k>r$, and note that for every arc $I\in\mathcal P_k$ if there exists $j$ for which $z_j^k\in I$, then for every $z\in I$ 
$$
\begin{cases}
\abs{\phi'\bb{z}} \in\bb{\delta_k^{1-\alpha+3\eta'},\delta_k^{1-\alpha-3\eta'}}\subset\bb{\bb{1-r_k}^{a+\eta},\bb{1-r_k}^{a-\eta}}\\
\abs{\phi'^{-i}\bb{z_j^k}}\in\bb{\delta_k^{\gamma+3\eta'},\delta_k^{\gamma-3\eta'}}\subset\bb{\bb{1-r_k}^{b+\eta},\bb{1-r_k}^{b-\eta}}
\end{cases},
$$
for $a=1-\frac1\alpha,\;\; b=\frac\gamma\alpha$, since $\log\phi'$ is a Bloch function, and 
$$
\delta_k^{1-\alpha\pm 3\eta'}=\bb{\delta_k^\alpha}^{1-\frac1\alpha\pm\frac{3\eta'}\alpha}=\bb{1-r_k}^{1-\frac1\alpha\pm\frac{3\eta'}\alpha}\;\;\text{ and }\;\;\delta_k^{\gamma\pm 3\eta'}=\bb{\delta_k^\alpha}^{\frac\gamma\alpha\pm\frac{3\eta'}\alpha}=\bb{1-r_k}^{\frac\gamma\alpha\pm\frac{3\eta'}\alpha}
$$
implying that for every $I\in\mathcal P_k$ containing $z_j^k$ for some $j$, $I\subset L_{a,b,\eta}(r_k)$. Next,
$$
\#\bset{I\in\mathcal P_k, \exists j \text{ s.t. } z_j^k\in I}=\#\bset{z_j^k}\gtrsim N(\delta_k,\alpha,\gamma,\eta').
$$
Overall,
\begin{align*}
\lambda_1\bb{L_{,b,\eta}(r_k)}\ge\bb{1-r_k}\#\bset{I\in \mathcal P_k, \exists j \text{ s.t. } z_j^k\in I}\gtrsim\bb{1-r_k}N(\delta_k,\alpha,\gamma,\eta')\\
\gtrsim\bb{1-r_k}\cdot\delta_k^{-\bb{f_\Omega(\alpha,\gamma)-\eta\cdot\alpha}}=\bb{1-r_k}^{1-\frac{f_\Omega(\alpha,\gamma)}\alpha+\eta},
\end{align*}
concluding the proof that
$$
d_\Omega(a,b)=\limsup_{\eta\searrow0}\;\limsup_{r\nearrow1^-}\frac{\log\bb{\lambda_1\bb{L_{a,b,\eta}(r_k)}}}{\log\bb{\frac1{1-r_k}}}+1\ge (1-a)f_\Omega\bb{\frac{1}{1-a},\frac{-b}{1-a}}.
$$

\paragraph*{Proof that $d_\Omega\le f_\Omega$}: Fix $\eps$ and let $\eta$ be small enough and $r$ close enough to $1$ so that $L_{a,b,\eta}(r)>\bb{\frac1{1-r}}^{d_\Omega(a,b)-\eps}$. Let $\bset{A_j(r)}$ be a collection of disjoint arcs of diameter between $\frac1C\cdot \bb{1-r}$ and $C\cdot \bb{1-r}$, for some uniform constant $C>0$, so that the collection $\bset{2A_j(r)}$ forms a covering for the set $L_{a,b,\eta}(r)$.  Because $\log\phi'$ is a Bloch function, such a cover always exists. Note that since the collection $\bset{2A_j(r)}$ covers $L_{a,b,\eta}(r)$ and it is minimal in the sense that $\bset{A_j(r)}$ are disjoint, then
 $$
 \#\bset{A_j(r)}= \#\bset{2A_j(r)}\ge \frac{L_{a,b,\eta}(r)}{4C(1-r)}>\bb{\frac1{1-r}}^{d_\Omega(a,b)-\eps}\cdot\frac1{4C}.
 $$

For every $j$ let $\xi_j$ denote the center of the arc $A_j$ and let $z_j:=\bb{1-r}\xi_j$. Then 
$$
\bb{\frac1{1-r}}^{a-\eta}< \abs{\phi'(z_j)} < \bb{\frac1{1-r}}^{a+\eta}\quad,\quad \bb{\frac1{1-r}}^{b-\eta}<\abs{\phi'(z_j)^{-i}}<\bb{\frac1{1-r}}^{b+\eta}.
$$
Following the Theorem \ref{thm:relations}, there exists a disk, $\tilde B_j=B(\zeta_j,\delta_j)$ containing $\phi(A_j)$ satisfying:
\begin{enumerate}[label=(\arabic*)]
\item  $\zeta_j\in\partial\Omega$
\item  There exists a constant $C=C(K_\Omega)$ so that $\frac{\delta_j}{\bb{1-r}^{1-a}}\in\bb{\frac{\bb{1-r}^\eta}{C},C\cdot\bb{1-r}^{-\eta}}$, i.e., $\delta_j$ is comparable with $(1-r)\abs{\phi'(z_j)}$.
\item $\tilde B_j$ satisfies \eqref{eq:derivative_measure} and \eqref{eq:derivative_rotation}.
\end{enumerate}
Let $\delta:= C\bb{1-r}^{1-a-\eta}$, then combining this with the inequalities above we see that
$$
\omega(\tilde B_j)\begin{cases}&<\frac{\delta_j}{\abs{\phi'(z_j)}}< C\bb{1-r}=\delta^{\frac1{1-a-\eta}}\cdot C^{\frac{a+\eta}{1-a-\eta}}\le \delta^{\frac1{1-a}}\\
&>\frac{\delta_j}{\log^2\bb{\frac1{\delta_j}}\abs{\phi'(z_j)}}>\frac{1-r}{C\log^2\bb{\frac1{\delta_j}}}>\frac{\delta^{\frac1{1-a-\eta}}}{\log^3\bb{\frac1{\delta_j}}}>\delta^{\frac{1}{1-a}+\sqrt\eta}
\end{cases}$$
and
$$
\rot(\tilde B_j)\in\bb{\frac1C\cdot  \abs{\phi'(z_j)^{-i}},C \abs{\phi'(z_j)^{-i}}}\subset\bb{\frac1C\bb{\frac1{1-r}}^{b-\eta},C\bb{\frac1{1-r}}^{b+\eta}}\subset\bb{\delta^{-\frac {b}{1-a}+\sqrt\eta},\delta^{-\frac{b}{1-a}-\sqrt\eta}},
$$
if $\delta$ is small enough as a function of $C,\; a$, and $b$.

The next problem is that our disks, $\tilde{B_j}$, have different radii, $\delta_j$, and they might not be disjoint. We therefore define $B_j:=B(\zeta_j,\delta)$, which is concentric with $\tilde B_j$.  The last part of the proof focuses on counting the number of disks in the collection $\bset{B_j}$ and showing for every $j$ the disk $B_j$ inherits the harmonic measure and the rotation of $\tilde B_j$ with controlled errors. 

We will first show that excluding an insignificant portion of the disks, we may assume that the initial collection, $\bset{\tilde B_j}$, is pairwise disjoint. Because harmonic measure of quasidisks is doubling (see Lemma \ref{lem:doubling_measure_quasidisk}), there exists a constant $M=M(K_\Omega)$ so that every disk, $\tilde B_j$, intersects at most 
\begin{align*}
&M\cdot C\bb{1-r}^{-\eta}\ge M\cdot\omega(\tilde B_j)\ge\omega(2\tilde B_j)\ge \underset{k, \tilde B_j\cap\tilde B_k\neq\emptyset}\sum\omega(\phi(A_k(r)))=\#\bset{k,\;\tilde B_j\cap\tilde B_k\neq\emptyset}\cdot \frac{1-r}C\\\
&\Rightarrow\#\bset{k,\;\tilde B_j\cap\tilde B_k\neq\emptyset}\le C\cdot M\cdot\frac{\omega(\tilde B_j)}{1-r}\lesssim \delta_j^{\frac1{1-a}-\eta}\cdot\delta^{-\frac1{1-a-\eta}}<\delta^{-\frac{2\eta}{(1-a)^2}}
\end{align*}
if $\eta$ is numerically small enough depending on $a$ and the constants of the quasidisk. By excluding at most $\delta^{\frac{2\eta}{(1-a)^2}}$ of the disks in the collection, we may assume without loss of generality that $\bset{\tilde B_j}$ are pairwise disjoint, and that $\#\bset{\tilde B_j}\ge \bb{\frac1{1-r}}^{d_\Omega(a,b)-2\eps}$ if $\eta$ is small enough as a function of $\eps$.

To bound the harmonic measure of $B_j$, note that by inclusion, 
$$
\omega(B_j)\ge\omega(\tilde B_j)\ge \delta^{\frac{1}{1-a}+\sqrt\eta}.
$$
For an upper bound however, we will use the fact that $\Omega$ is a quasidisk, or, more importantly, that its harmonic measure is doubling. Following property (2) of $\tilde B_j$, we have that since $\delta= C\bb{1-r}^{1-a-\eta}$, then$$
\delta>\delta_j>\frac{\bb{1-r}^{1-a+\eta}}C=\bb{\frac\delta C}^{\frac{1-a+\eta}{1-a-\eta}}>\delta^{1+\frac{4\eta}{1-a}}=\delta\cdot\delta^{\frac{4\eta}{1-a}}
$$
if $\eta$ is small enough as a function of the uniform constant $C$, and $(1-a)$. In particular,
$$
2^{-m}B_j\subseteq \tilde B_j\subseteq B_j,
$$
where $m:=\left\lceil{\log_2\bb{\delta^{-\frac{4\eta}{1-a}}}}\right\rceil$. Since $B_j$ is centred at $\partial\Omega$, following Lemma \ref{lem:doubling_measure_quasidisk}
$$
\omega(B_j)\le C\omega\bb{\frac12 B_j}\le\cdots\le C^m\omega\bb{2^{-m}B_j}\le C^m\omega(\tilde B_j)\le C^m\cdot \delta^{\frac1{1-a}}=C^{\left\lceil{\log_2\bb{\delta^{-\frac{4\eta}{1-a}}}}\right\rceil}\delta^{\frac1{1-a}}\le \delta^{\frac1{1-a}-\sqrt\eta},
$$
if $\eta$ is small enough as a function of the uniform constant $C$, and $(1-a)$.

To show that the rotation of $B_j$ does not change much, we will use Lemma \ref{lem:rot_quasidisk_radius}, with $B_1=B_j$ and $B_2=\tilde B_j$. Note that the radius of $B_j$ is $\delta$ and the radius of $\tilde B_j$ is at least $\delta\cdot\delta^{\frac{4\eta}{1-a}}$, i.e., 
\begin{align*}
\abs{\log\bb{\rot(B_j)}-\log\bb{\rot(\tilde B_j)}}\le 120  K_\Omega\cdot\frac{\log\bb{\frac1{\delta_2}}}{\log\bb{\frac1{\delta_1}}}\le 120  K_\Omega\cdot\frac{\log\bb{\frac1{\delta\cdot\delta^{\frac{4\eta}{1-a}}}}}{\log\bb{\frac1{\delta}}}\\
=120  K_\Omega\cdot\frac{\log\bb{\frac1{\delta}\bb{1+\frac{4\eta}{1-a}}}}{\log\bb{\frac1{\delta}}}\le 2\sqrt\eta\log\bb{\frac1\delta}.
\end{align*}
Then, combined with the bounds on the rotation of $\tilde B_j$, we obtain that
$$
\rot(B_j)\begin{cases}
		\le \rot(\tilde B_j)\cdot\exp\bb{\abs{\log\bb{\rot(B_j)}-\log\bb{\rot(\tilde B_j)}}}\le \delta^{-\frac{b}{1-a}-2\sqrt\eta}\\
		\ge \rot(\tilde B_j)\cdot\exp\bb{-\abs{\log\bb{\rot(B_j)}-\log\bb{\rot(\tilde B_j)}}}\ge \delta^{-\frac{b}{1-a}+2\sqrt\eta}
	\end{cases}.
$$

Next, while the disks in the collection $\bset{\tilde B_j}$ are disjoint, $B_j\supseteq\tilde B_j$ which means that the collection $\bset{B_j}$ might not be a-priori disjoint. Fix $j$ and note that if $B_j\cap B_k\neq\emptyset$ then $\abs{z_k-z_j}<3\delta$, i.e., $\tilde B_k\subseteq B_k\subset 3B_j$.  This implies that since the collection $\bset{\tilde B_k}$ is disjoint,
\begin{align*}
9\pi\cdot\delta^2&=\lambda_2\bb{3B_j}\ge\lambda_2\bb{\underset{B_k\cap B_j\neq\emptyset}\bigcup \tilde B_k}=\#\bset{k, B_j\cap B_k\neq\emptyset}\cdot\underset k\inf\; \lambda_2\bb{\tilde B_k}\ge\#\bset{k, B_j\cap B_k\neq\emptyset}\cdot\pi\delta^{2\bb{1+\frac{4\eta}{1-a}}}\\
&\Rightarrow \#\bset{k, B_j\cap B_k\neq\emptyset}\le\frac{9\pi\cdot\delta^2}{\pi\delta^{2\bb{1+\frac{4\eta}{1-a}}}}=9\delta^{-\frac{8\eta}{1-a}}.
\end{align*}
As before, by excluding at most $9\delta^{\frac{8\eta}{1-a}}$ of the disks in the collection, we may assume without loss of generality that $\bset{ B_j}$ are pairwise disjoint as well. Overall, the collection $\bset{B_j}$ is a collection of pairwise disjoint disks centered at $\partial\Omega$ and satisfying \eqref{eq:N_measure} and \eqref{eq:N_rotation}, implying that
\begin{align*}
f_\Omega\bb{\frac{1}{1-a},\frac{-b}{1-a}}= \limit \eta 0 \limitsup \delta 0 \frac{\log N\bb{\delta,\frac{1}{1-a},\frac{-b}{1-a},2\sqrt\eta}}{\log\bb{\frac1\delta}}\ge \limit \eta 0 \limitsup \delta 0 \frac{\log\bb{\#\bset{B_j}}}{\log\bb{\frac1\delta}}\\
\ge\limit \eta 0 \limitsup \delta 0 \frac{ \log\bb{\frac{\#\bset{\tilde B_j}}{9\delta^{\frac{8\eta}{1-a}}}}}{\log\bb{\frac1\delta}}\ge \limit \eta 0 \limitsup \delta 0\frac{ \log\bb{\frac{\bb{\frac1{\delta}}^{\frac{d_\Omega(a,b)-2\eps}{1-a-2\eta}}}{9\delta^{\frac{8\eta}{1-a}}}}}{\log\bb{\frac1\delta}} \ge \frac{d_\Omega(a,b)-3\eps}{1-a}
\end{align*}
since $\eps>0$ was generic, taking $\eps\searrow 0$ concludes the proof.

\paragraph*{The Proof $\displaystyle{\min_{\sigma,\ \sigma'\in\bset{-,+}}\bset{d_\Omega^{\sigma\sigma'}}}$ and $\displaystyle{\min_{\sigma,\ \sigma'\in\bset{-,+}}\bset{f_\Omega^{\sigma\sigma'}}}$ satisfy Relation \eqref{eq:distdim}:}~\\
The relation follows from the fact that at the growth points of the corresponding functions, they are equal to $d_\Omega$ and $f_\Omega$, respectively. Moreover, the functions can be restored from their values at those points.

%% file: sections/spectra.tex
In this section, we present the Minkowski Crosscut Spectrum. This novel kind of spectrum relates counting crosscuts with certain harmonic measure and certain rotation with the Minkowski Distortion Spectrum.

 In the context of the Main Lemma, Lemma \ref{lem:main}, it is instructive to consider the implications of enumerating crosscuts with disjoint supports in place of disjoint disks. This alternative approach prompts a deeper examination of the spectral characteristics associated with crosscuts. Accordingly, we extend the discussion of spectral frameworks by analyzing the spectrum of crosscuts and elucidating its connection to the Miskowski Distortion Mixed Spectrum through an appropriate collection of such crosscuts.

For the reader's convenience we repeat the definition:

\begin{defn*}
For every $r\in\bb{0,1}$ and $a,b>0$ we define by $\Gamma\bb{a,b,r}$ to be a maximal collection of crosscuts, $C$, with disjoint supports satisfying
\begin{enumerate}
\item $\abs{\frac{\omega(\nu_C)}{1-r}-1}\le\frac{1-r}{\log\log\log\bb{\frac1{1-r}}}$.
\item $\diam(\nu_C)\ge\bb{1-r}^{1-a}$.
\item $\rot(C)\ge\bb{1-r}^{-b}$.
\end{enumerate}
We define the \emph{Minkowski Crosscut Spectrum} by
$$
d_\Omega^{cc}(a,b)=\underset{a'\searrow a\atop {b'\searrow b}}\limsup\;\underset{r\nearrow 1}\limsup\;\frac{\log\bb{\# \Gamma\bb{a',b',r}}}{\log\bb{\frac1{1-r}}}.
$$
\end{defn*}

In a sense, if $\alpha=\frac1{1-a}$, $\gamma=\frac{-b}{1-a}$, and $1-r\sim\eps^\alpha$ then every crosscut, $C$, in this collection satisfy $\omega(\nu_C)\sim\eps^\alpha$, $\diam(\nu_C)\ge\eps$, and $\rot(C)\ge\eps^\gamma$, where $\nu_C$ denotes the support of the crosscut, $C$.



The first Lemma in this section shows that there is some correspondence between the Minkowski Crosscut Spectrum, $d_\Omega^{cc}$, and the classical Minkowski Mixed Distortion Spectrum, $d_\Omega$. Recall that by the symmetry argument we can always assume that $b\ge 0$.
\begin{lem*}[\ref{lem:curves_bnd_d}]
Let $a>0, b\in\R$.
\begin{enumerate}
\item  For every simply connected domain, $\Omega$, $d_\Omega(a,b)\le d_\Omega^{++}(a,b)\le d_\Omega^{cc}(a,b)$.
\item For every quasidisk, $\Omega$, $d_\Omega(a,b)= d_\Omega^{cc}(a,b)$.
\end{enumerate}
\end{lem*}
\begin{rmk} While the Minkowski Crosscut Spectrum majorates the Minkowski Distortion Mixed Spectrum, it still does not \underline{always} majorate the Minkowski Dimension Mixed Spectrum, as \cite[Theorem 2.3]{AdiBin2} shows.\end{rmk}
{\bf Proof of \eqref{item:ineq_ds}:} Let $\eta>0$, fix $r$ close enough to 1, and let
$$
L^{++}_{a,b,\eta}(r)=\bset{\zeta\in\partial\D, \log\abs{\varphi'(r\zeta)}>\bb{a-\eta}\log\bb{\frac1{1-r}},\; \argu\bb{\phi'(r\zeta)}>\bb{b-\eta}\log\bb{\frac1{1-r}}}.
$$
Let $\Gamma=\bset{A_j}_{j=1}^M$ denote the minimal collection of disjoint arcs in $\partial\D$ satisfying that $\lambda_1(A_j)=1-r$ while $\bunion j 1 M 2A_j\supset L^{++}_{a,b,\eta}(r)$. For every $j$ there exists $\zeta'\in 2A_j\cap L^{++}_{a,b,\eta}(r)\neq\emptyset$, as $\bset{A_j}$ is minimal. Let $a_j$ denote the centre of the arc $A_j$. The function $\phi$ is conformal, making $\log\phi'$ a Bloch function. Therefore, 
$$
\max\bset{\abs{\log\abs{\phi'(r\cdot \zeta)}-\log\abs{\phi'(r\cdot a_j)}},\abs{\argu\bb{\phi'(r\cdot \zeta)}-\argu\bb{\phi'(r\cdot a_j)}}}\le \abs{\log\bb{\phi'(r\cdot \zeta)}-\log\bb{\phi'(r\cdot a_j)}}\le 6,
$$
and for every $\zeta\in 2A_j\cap L^{++}_{a,b,\eta}(r)$
\begin{align*}
&\log\abs{\phi'(r\cdot a_j)}\ge\log\abs{\phi'(r\cdot \zeta)}-6\ge \bb{a-\eta}\log\bb{\frac1{1-r}}-6\ge\bb{a-2\eta}\log\bb{\frac1{1-r}}\\
&\argu\bb{\phi'(r\cdot a_j)}\ge \argu\bb{\phi'(r\cdot \zeta)}-6\ge \bb{b-\eta}\log\bb{\frac1{1-r}}-6\ge\bb{b-2\eta}\log\bb{\frac1{1-r}}.
\end{align*}
Let $C_j$ be the crosscut of the arc $A_j$. Following Koebe's distortion theorem combined with the first part of Observation \ref{obs:distortion_vs_curve},
$$
\diam(\nu_{C_j})\gtrsim\abs{\phi'(a_j)}\bb{1-\abs{a_j}}\ge \bb{\frac1{1-r}}^{a-\eta-\frac{C}{\log\bb{\frac1{1-r}}}-1}=\bb{1-r}^{1-a+\eta+\frac{C}{\log\bb{\frac1{1-r}}}}>\bb{1-r}^{1-a+2\eta},
$$
if $r$ is close enough to 1, depending on $\eta$, where $C$ is some absolute constant. In addition, by the way the arcs $A_j$ were defined
$$
\rot(C_j)\ge \frac{\abs{\phi'^{-i}\bb{r\cdot a_j}}}{C_\Omega\bb{1-r}^{c(r)}}\ge\bb{1-r}^{-b+\eta+c(r)}\ge\bb{1-r}^{-b+2\eta}\quad;\quad \omega(\nu_{C_j}) \begin{cases}\ge \lambda_1(A_j)=1-r \\
																										\le 2\lambda_1(A_j)=2(1-r)
																							\end{cases}
$$
for some $c(r)\searrow 0$ as $r\nearrow 1^-$, following the main lemma, Lemma \ref{lem:main} ($C_\Omega$ a-priori refers to two different constants, that by taking the maximum one can be unified).

Define $a'=a-2\eta,\; b'=b-2\eta$, then $C_j\in \Gamma\bb{a',b',r}$. In particular,
$$
M(\eta,r)=\text{number of curves }C_j\le\# \Gamma\bb{a',b',r},
$$
while $a'\rightarrow a$ and $b'\rightarrow b$ as $r\rightarrow 1$ and $\eta\searrow 0$. We conclude that
\begin{align*}
d_\Omega^{++}\bb{a,b}&=\limsup_{\eta\searrow 0}\;\limsup_{r\nearrow1^-}\frac{\log\bb{\lambda_1\bb{L^{++'}_{a,b,\eta}(r)}}}{\log\bb{\frac1{1-r}}}+1\le \limsup_{\eta\searrow 0}\;\limsup_{r\nearrow1^-}\frac{\log\bb{2(1-r)\cdot M(\eta,r)}}{\log\bb{\frac1{1-r}}}+1\\
&\le \underset{a'\searrow a\atop b'\searrow b}\limsup\;\underset{r\nearrow 1}\limsup\;\frac{\log\bb{\# \Gamma\bb{a',b',r}}}{\log\bb{\frac1{1-r}}}\le d_\Omega^{cc}(a,b),
\end{align*}
concluding the proof of \eqref{item:ineq_ds}. Note that if $\Omega$ is a quasidisk, all the constants in the estimates above depend only on the dilatation, $K_\Omega$.

{\bf Proof of \eqref{item:quasidisk_curve}:} In light of \eqref{item:ineq_ds}, we only need to show that $d_\Omega^{cc}(a,b)\le d_\Omega(a,b)$. 

Let $C\in\Gamma(a,b,r)$ be a crosscut, then it is supported on a curve, $\nu_C$, and $I_C:=\phi\inv(\nu_C)\subset\T$ is an arc. Note that, by definition, $C$ is a crosscut of $I_C$. We denote by $\widetilde z_C$ the point representing $I_C$, and note that, since $\Omega$ is a quasidisk, Observation \ref{obs:distortion_vs_curve} shows that the estimates done in \eqref{item:ineq_ds} are equalities and not inequalities, i.e., 
$$
\bb{1-r}^{-a}\le\diam(\nu_C)\dsim{K_\Omega}\abs{\phi'(\widetilde z_C)}\bb{1-\abs{\widetilde z_C}}\Rightarrow \abs{\phi'(\widetilde z_C)}\dmore{K_\Omega} \bb{1-\abs{\widetilde z_C}}^{-a}.
$$
In addition,
$$
\abs{\frac{1-\abs{\widetilde z_C}}{1-r}-1}=\abs{\frac{\omega(\nu_C)}{1-r}-1}<\frac{1-r}{\log\log\log\bb{\frac1{1-r}}}\quad\quad\text{ and }\quad\quad\bb{1-r}^{-b}\le \rot(C)\sim \abs{\phi'^{-i}(\widetilde z_C)}.
$$
Let $z_C=r\cdot\frac{\widetilde z_C}{\abs{\widetilde z_C}}$. Then, since $\log\phi'$ is a Bloch function,
\begin{align*}
&\abs{\log\bb{\phi'(z_C)}-\log\bb{\phi'(\widetilde z_C)}}\le 6arctanh\bb{\rho(z_C,\widetilde z_C)}\le 12\frac{\abs{z_C-\widetilde z_C}}{\abs{1-z_C\cdot \overline{\widetilde z_C}}}\lesssim\frac{1}{\log\log\log\bb{\frac1{1-r}}}\\
&\Rightarrow \begin{cases}
				\abs{\phi'(z_C)}\ge\abs{\phi'(\widetilde z_C)}\bb{1-\frac{Const}{\log\log\log\bb{\frac1{1-r}}}}\ge \bb{1-r}^{-a+\eta}\\
				\abs{\phi'^{-i}( z_C)}\ge\abs{\phi'^{-i}(\widetilde z_C)}\bb{1-\frac{Const}{\log\log\log\bb{\frac1{1-r}}}}\ge \bb{1-r}^{-b+\eta}
			\end{cases}
\end{align*}
for some $\eta>0$ small. Note that  since crosscuts in $\Gamma(a,b,r)$ have disjoint supports, each crosscut gives rise to a unique arc of length $(1-r)$ contained in $L^{++'}_{a,b,\eta}(r)$. We conclude that
\begin{align*}
d_\Omega^{++}\bb{a,b}&=\limsup_{\eta\searrow 0}\;\limsup_{r\nearrow1^-}\frac{\log\bb{\lambda_1\bb{L^{++'}_{a,b,\eta}(r)}}}{\log\bb{\frac1{1-r}}}+1\ge\limsup_{r\nearrow1^-}\frac{\log\bb{\#\Gamma\bb{a,b,r}}}{\log\bb{\frac1{1-r}}}=d_\Omega^{cc}(a,b).
\end{align*}
\begin{align*}
\quad\quad\quad\quad\quad\quad\quad\quad\quad\quad\quad\quad\quad\quad\quad\quad\quad\quad\quad\quad\quad\quad\quad\quad\quad\quad\quad\quad\quad\quad\quad\quad\quad\quad\quad\quad\quad\quad\quad\quad\quad\quad\quad\quad\quad\quad\quad\quad\qed
\end{align*}

We are not using the Minkowski Crosscut Spectrum for $a<0$. Instead, we use the following lemma to show that for $a<0$ there is no need to use the Minkowski Crosscut Spectrum, as, in this case, the Minkowski Dimension Mixed Spectrum, $f^{++}$, majorates the Minkowski Distortion Mixed Spectrum.
\begin{lem*}[\ref{lem:neg_a}]
If $a<0$, then for every simply connected domain, $\Omega$,
$$
d_\Omega(a,b)\le d_\Omega^{++}(a,b)\le \limitinf\eta{0^+} f_\Omega^{+\pm}\bb{\frac{1}{1-a}+\eta,\frac{-b}{1-a}+\eta}.
$$
\end{lem*}
\begin{proof}
While we prove the upper bound for $f^{++}$ the same proof with very mild changes shows the same for $f^{+-}$. We will first prove the statement for growth points. The proof of this part is composed of two steps; defining the disks and proving they satisfy \eqref{eq:N_measure} and \eqref{eq:N_rotation}.

Fix a growth point, $(a,b)$, and $\eps_1$ and let $\eta_0>0$ be small enough so that for every $r$ large enough $r>r(\eta_0)$
$$
\lambda_1\bb{L_{a,b,\eta_0}(r)}>(1-r)^{1+\eps_1-d_\Omega(a,b)}.
$$
Since $(a,b)$ is a growth point, 
$$
d_\Omega^{++}(a-\eta_0,b-\eta_0)>d_\Omega^{++}(a+\eta_0,b+\eta_0).
$$
We define 
$$
\eps:=\min\bset{\eps_1,\frac13\bb{d_\Omega^{++}(a-\eta_0,b-\eta_0)-d_\Omega^{++}(a+\eta_0,b+\eta_0)}}.
$$
Fix $\eta<\eta_0$ be small enough and $r<r(\eta,\eta_0)$ close enough to 1 so that $\bb{1-r}^\eps<\frac12$ and 
$$
\lambda_1\bb{L^{++}_{a-\eta_0,b-\eta_0,\eta}(r)}>(1-r)^{1+\eps-M}>(1-r)^{1+2\eps-M}>\lambda_1\bb{L^{++}_{a+\eta_0,b+\eta_0,\eta}(r)}.
$$
for some $M-d_\Omega^{++}(a,b)\in\bb{-\eps_1,\eps_1}$.

Define $(a',b')=(a-\eta_0,b-\eta_0)$, partition $\T$ into $n:=\left\lceil\frac1{1-r}\right\rceil$ arcs, and let $\bset{I_j}, I_j\subset\partial\D$ be the minimal collection of arcs in the partition with $L^{++}_{a',b',\eta}(r)\subset \bigcup_j 2I_j$. In particular, for every $j$ we have $I_j\cap L^{++}_{a',b',\eta}(r)\neq\emptyset$, which implies that for every $\zeta\in I_j$
$$
\abs{\phi'(r\cdot\zeta)}\ge C\bb{1-r}^{\abs {a'}+\eta}\quad\text{and}\quad \abs{\phi'^{-i}(r\zeta)}\ge \frac1C(1-r)^{-b'+\eta}
$$
for some uniform constant $C$, and, in particular, this holds for the center of $I_j$, denoted $\zeta_j$. 

Let $J_0$ denote the sub-collection of indices $j$ so that
$$
\frac1C(1-r)^{\abs{a'}+\eta}\le\abs{\phi'(r\cdot\zeta_j)}< C\bb{1-r}^{\abs {a'}-2\eta_0+\eta}\quad,\quad\frac1C(1-r)^{ -b'+\eta}<\abs{\phi'^{-i}(r\cdot\zeta_j)}\le C\bb{1-r}^{- b'+2\eta_0-\eta}.
$$
i.e., these are indices so that elements in $I_j$ are in $L^{++}_{a',b',\eta}(r)=L^{++}_{a-\eta_0,b-\eta_0,\eta}(r)$ but not in $L^{++}_{a'+2\eta_0,b+2\eta_0,\eta}(r)=L^{++}_{a+\eta_0,b+\eta_0,\eta}(r)$. Note that  since $r$ is close enough to 1,
\begin{align*}
\# J_0&=\frac{\lambda_1\bb{\underset{j\in J}\uplus I_j}}{\frac1n}\ge n\cdot\lambda_1\bb{L^{++}_{a',b',\eta}(r)\setminus L^{++}_{a'+2\eta_0,b',\eta}(r)}=n\bb{\lambda_1\bb{L^{++}_{a-\eta_0,b-\eta_0,\eta}(r)}-\lambda_1\bb{L^{++}_{a+\eta_0,b+\eta_0,\eta}(r)}}\\
&\ge n\bb{(1-r)^{1+\eps-M}-(1-r)^{1+2\eps-M}}=n(1-r)^{1+\eps-M}\bb{1-(1-r)^{\eps}}\ge\frac14 (1-r)^{-M}>(1-r)^{2\eps_1-d_{\Omega}^{++}(a,b)}
\end{align*}
by the way all the parameters were defined and chosen.

Let $j\neq k$, and assume without loss of generality that $\abs{\phi'(z_j)}\ge\abs{\phi'(z_k)}$, then, following \cite[Cor 1.5]{pombook},
\begin{equation}\label{eq:dist}
\abs{\phi(z_j)-\phi(z_k)}\ge \bb{1-\abs{z_j}^2}\abs{\phi'(z_j)}\abs{\frac{z_j-z_k}{1-z_j\overline{z_k}}}\ge\abs{\phi'(z_j)}\cdot \abs{z_j-z_k}\gtrsim\abs{z_j-z_k}\cdot (1-r)^{\abs{a'}+\eta}.
\end{equation}
We define the collection of indices, $J$, by excluding at most $(1-r)^{3\eta_0}$ of the indices in $J$ to guarantee that whenever $j,k\in J$ we have $\abs{\phi(z_j)-\phi(z_k)}>(1-r)^{1+\abs{a'}-3\eta_0}$.

For every $j\in J$, let $C_j$ be a crosscut of the interval $I_j$. Note that
$$
\diam(C_j)\sim \diam(\phi(I_j))+2\dist(\phi(I_j),\partial\Omega)\in\bb{\frac1{C_\Omega}(1-r)^{1+\abs{a'}+\eta},C_\Omega(1-r)^{1+\abs{a'}-2\eta_0+\eta}}.
$$
Define $z_j=r\cdot\zeta_j$, and let $\xi_j\in\partial\Omega$ be the closest point to $\phi(z_j)$. We may assume without loss of generality that $\phi(z_j)\in C_j$ and that $\xi_j$ is in the support of the crosscut $C_j$. Following the main lemma, Lemma \ref{lem:main},
$$
\abs{\phi'\bb{z_j}} \dsim \Omega \frac{\diam(C_j)}{\lambda_1(I_j)}\quad\quad;\quad\quad \frac{\log\abs{\phi'^{-i}\bb{z_j}}}{\log\rot(C_j)}\dsim \Omega 1.
$$ 
For every $j$ we define $\delta_j:=\dist(\xi_j,C_j)+\diam(C_j)$, then
$$
\delta_j 	\begin{cases}
			\le \diam(C_j)+\abs{\phi(z_j)-\xi_j}\le C_\Omega\bb{1-r}^{1+\abs a'-2\eta_0+\eta}\\
			\ge\diam(C_j)\ge \lambda_1(I_j)\abs{\phi'(z_j)}\ge \frac1{C_\Omega}(1-r)^{1+\abs{a'}+\eta}
		\end{cases}
$$
since for every $j\in J$ 
$$
\dist(C_j,\xi_j)\le \abs{\xi_j-\phi(z_j)}\le \bb{1-r^2}\abs{\phi'(z_j)}<C\bb{1-r}^{1+\abs a'-2\eta_0+\eta}.
$$
The constant $C_\Omega$ might be different but it still depends on the domain alone. 

Let $\delta=\underset{j\in J}\min\; \delta_j$, then for every $j$
$$
1\ge\frac{\delta}{\delta_j}\ge \frac{\frac1{C_\Omega}(1-r)^{1+\abs{a'}+\eta}}{C_\Omega\bb{1-r}^{1+\abs a'-2\eta_0+\eta}}=C_\Omega^{-2}(1-r)^{2\eta_0}.
$$
First we will show that the disks, $B_j=B(\xi_j,\delta)$ are disjoint; let $j\neq k$, following \eqref{eq:dist},
\begin{align*}
\abs{\xi_j-\xi_k}&\ge\abs{\phi(z_j)-\phi(z_k)}-\abs{\xi_j-\phi(z_j)}-\abs{\xi_k-\phi(z_k)}\\
&\ge (1-r)^{1+\abs{a'}-3\eta_0}-2C_\Omega\bb{1-r}^{1+\abs a'-2\eta_0+\eta}>2\max\bset{\delta_j,\delta_k}>2\delta,
\end{align*}
i.e., the disks $\bset{B_j}$ are pairwise disjoint.

Next, we will bound the harmonic measure of $B_j$ from below; let $\Omega_r=\phi(r\D)$, then $\phi(z_j)\in\partial\Omega_r$ and moreover for every $\zeta\in r\cdot I_j$,
$$
\dist(\phi(\zeta),B_j\cap\phi(I_j))\le \abs{\xi_j-\phi(z_j)}+\diam(\phi( I_j))\lesssim\abs{\phi'(z_j)}\diam(I_j)\le C\cdot (1-r)^{1+\abs{a'}}.
$$
In particular, Buerling's projection theorem yields $\omega(\zeta,B_j\cap\phi(I_j);\Omega)\gtrsim (1-r)^{2\eta_0} $ with uniform constants. Since the function $z\mapsto\omega(z,B_j\cap\phi(I_j);\Omega)$ is harmonic in $\Omega_r$ we obtain
\begin{align*}
\omega(B_j)\ge\omega(z_0,B_j\cap\phi(I_j);\Omega)=\integrate{\partial\Omega_r}{}{\omega(\zeta,B_j\cap\phi(I_j);\Omega)}\omega(z_0.\zeta;\Omega_r)\gtrsim  (1-r)^{2\eta_0} \omega(z_0,\phi(rI_j),\Omega_r)\\
= (1-r)^{2\eta_0} \lambda_1(r\cdot I_j)\ge (1-r)^{1+3\eta_0} \ge \delta^{\frac{1+3\eta_0}{1+\abs{a'}-2\eta_0}+\eta}.
\end{align*}
We conclude that \eqref{eq:N_measure}(+) holds as long as $r>r(\eta)$ is large enough enough.

To prove $B_j$ satisfies \eqref{eq:N_rotation}, let $B_1=B(w,\diam(C_j))$ be any disk centered at $C_j$. There exists a disk, $B_0$, of radius $\diam(C_j)$ so that $\phi(z_j)\in\partial B_0$. Since $\phi(z_j)\in\partial B_0$, it can be chosen to contain at least half of the crosscut $C_j$, implying that $\frac{\diam(\partial\Omega\cap B_0)}{\diam(B_0)}\ge \frac{\diam(C_j)\bb{2-\frac12}}{2\diam(C_j)}>\frac14$, i.e., the rotation of $B_0$ is well defined, and $B_1\cap B_0\neq\emptyset$. Using Lemma \ref{lem:rotation_err_in}, and the definition of rotation of a crosscut, $\rot(C_j)\sim\rot(B_1)\sim \rot(B_0)$ with uniform constants.

Let $\nu\subset\Omega_{B_0}$ be a curve connecting $z_0$ with $\phi(z_j)\in \partial B_0$, where  $\Omega_{B_0}$ is the connected component of $\Omega\setminus B_0$ containing $z_0$. We may assume without loss of generality that $\nu\cap B_j=\emptyset$, since $\bb{B_j\cap \nu_{C_j}}\cap B_0=\emptyset$. We define the curve $\nu'=\nu+\bb{\sbb{\phi(z_j),\xi_j}\setminus B_j}$, then $\nu'\subset\Omega_{B_j}$ since $\nu\subset \Omega_{B_j}$ and $\sbb{\phi(z_j),\xi_j}\setminus B_j\subset\Omega_{B_j}$, as $\xi_j$ is the closest point to $\phi(z_j)$ in $\partial\Omega$ implying that $\sbb{\phi(z_j),\xi_j}\subset B\bb{\phi(z_j),\abs{\xi_j-\phi(z_j)}}\subset\Omega$. Combining everything together and using Lemma \ref{lem:rotation_integral}, we see that
\begin{align*}
\abs{\log\rot(C_j)-\Im\int_\nu\frac 1{z-\xi_j}dz}\le \abs{\log\rot(C_j)-\log\rot(B_0)}+\abs{\log\rot(B_0)-\Im\int_\nu\frac 1{z-\xi_j}dz}+C\\
\le\abs{\log\rot(B_0)-\Im\int_{\nu}\frac 1{z-w_0}dz}+\abs{\Im\int_{\nu}\frac 1{z-w_0}dz-\Im\int_\nu\frac 1{z-\xi_j}dz}+C\le C,
\end{align*}
as the winding index of $\nu$ around $\xi_j$ and the winding index of $\nu$ around $w_0$ is the same, i.e., 
$$
\frac1{2\pi i}\underset{ \nu}\int \frac1{z-w_0}dz=\frac1{2\pi i}  \underset{ \nu}\int \frac1{z-\xi_j}dz.
$$
Note that the constant $C$ is a changing from line to line, however it is uniform. Finally, since $\sbb{\phi(z_j),\xi_j}\setminus B_j$ is a straight line, then the rotation along it changes by at most $\pi$ implying that
\begin{align*}
\abs{\log\rot(C_j)-\log\rot(B_j)}\le \abs{\log\rot(C_j)-\Im\int_\nu\frac 1{z-\xi_j}dz}+\abs{\Im\int_{\nu'}\frac 1{z-\xi_j}dz-\log\rot(B_j)}\\
+\abs{\Im\int_\nu\frac 1{z-\xi_j}dz-\Im\int_{\nu'}\frac 1{z-\xi_j}dz}\lesssim 1
\end{align*}
with uniform constants.

We conclude that the disks $\bset{B_j}$ are pairwise disjoint, centered at $\partial\Omega$, and satisfy both \eqref{eq:N_measure} and \eqref{eq:N_rotation}, i.e., $N^{++}\bb{\delta,\frac{1+3\eta_0}{1+\abs{a'}-2\eta_0},\frac{-b}{1+\abs{a'}-2\eta_0},2\eta}\ge\#\bset{B_j}$, implying that
$$
 f_\Omega^{++}\bb{\frac{1+3\eta_0}{1-a-\eta_0},\frac{-b}{1-a-\eta_0}}\ge d_{\Omega}^{++}(a,b)-2\eps_1-3\eta_0,
 $$
taking $r\nearrow 1$ (which imposes $\delta\searrow 0$) and $\eta\searrow 0$. If we then take $\eta_0\searrow 0$ and $\eps_1\searrow 0$ we obtain
\begin{align*}
\limitinf\eta{0^+} f_\Omega^{++}\bb{\frac{1}{1-a}+\eta,\frac{-b}{1-a}+\eta}\ge  d_{\Omega}(a,b)^{++},
\end{align*}
since $f^{++}$ is monotone increasing, and upper semi-continuous.

Finally, if $(a,b)$ is not a growth point of $d$, let $(a',b')$ be a growth point satisfying $a'>a, b'>b$. Then
$$
d_\Omega(a,b)^{++}=d_\Omega^{++}(a',b')\le f_\Omega^{++}\bb{\frac1{1-a'},\frac{-b'}{1-a'}}\le  f_\Omega^{++}\bb{\frac1{1-a},\frac{-b}{1-a}}
$$
since $f^{++}$ is monotone increasing.
For bounds on $f^{+-}$ we take $a'>a, b'<b$ and use monotonicity in two steps, i.e., 
$$
f_\Omega^{+-}\bb{\frac1{1-a'},\frac{-b'}{1-a'}}\le  f_\Omega^{+-}\bb{\frac1{1-a'},\frac{-b}{1-a'}}\le f_\Omega^{+-}\bb{\frac1{1-a},\frac{-b}{1-a}}
$$
concluding the proof.
\end{proof}

%% file: sections/repellers.tex
\subsection{Thermodynamic formalism}\label{subsec:formalism}
For the reader's convenience we repeat the definition of Jordan Repellers.
\begin{defn*}
Let $\Omega$ be a simply connected domain, and $\partial$ be a sub-arc of its boundary, $\partial\Omega$. We call $\partial$ a \emph{Jordan Repeller} if there exists a partition of $\partial$ into a finite number of non-intersecting sub-arcs, denoted $\partial_1,\partial_2,\cdots,\partial_N$ (a Markov partition), and a piecewise univalent (though not necessary injective) map $F$, such that for every sub-arc in the partition, $\partial_j$, there exists a neighborhood $U_j$ containing $\partial_j$ and a univalent map $F_j:\, U_j\rightarrow\C$ satisfying:
\begin{enumerate}[label=(\arabic*)]
\item $F|_{\partial_j}=F_j$.
\item (Geometry invariance) $F_j\bb{\Omega\cap U_j}\subset\Omega$.
\item (Boundary invariance) $F_j\bb{\partial\Omega\cap U_j}\subset\partial\Omega$.
\item (Markov property) The image of each $\partial_j$ is a finite union arcs $\partial_k$, i.e., $F(\partial_j)=\underset{k\in\mathcal A_j}\bigcup\partial_k$, where $\mathcal A_j\subseteq\bset{1,2,\cdots,N}$.
\item (Expanding) There exists $n_0$ so that $\underset{z\in\partial}\inf\;\abs{F^{(n_0)} (z)} >1$, where $F^{(n_0)}$ denotes the $n_0$'th derivative of $F$.
\item (Mixing) For any open set $W$ satisfying $W\cap\partial\neq\emptyset$, there exists $n$ such that $F^{\circ n}\bb{W\cap\partial}=\partial$, i.e., the $n$'th iteration of $F$ maps $W\cap\partial$ to cover $\partial$.
\end{enumerate}
Note that $F$ can be a multifold only near the end-points of the arcs $\partial_j$.
\end{defn*}
\subsection{On Symbolic Representations of Repellers}
The function $F$ is locally invertible, i.e., on the domains $\bset{U_j}_{j=1}^N$, the function $F|_{U_j}$ is invertible. Following Mauldin - Urbanski in \cite{MaulBanski1996}, if we denote by $f_j=\bb{F|_{U_j}}\inv$, then $f_j$ is contracting (as $F$ was expanding as a repeller). We define for every $j$ the set $V_j=\phi(U_j\cap\D)$ and note that $\partial_j=\partial V_j\cap\partial\Omega$. For every finite word, $X=(x_1,x_2,\cdots,x_n)$, over the alphabet $\Sigma=\Sigma_F=\bset{1,2,\cdots,N}$ we define the map 
$$
f_X:=f_{x_1}\circ f_{x_2}\circ\cdots\circ f_{x_n},
$$ 
and for every (in)finite word we denote by $X|_n=(x_1,\cdots,x_n)$.

Since $\bset{f_{x_j}}$ are all contracting, and the set $\partial\Omega$ is compact, for every infinite word, $X$, the set
$$
\underset{n\in\N}\bigcap f_{X|_n}(\partial\Omega)
$$
is a singleton. We define the map $\pi:\Sigma^\N\rightarrow\partial\Omega$, matching infinite words and points on the boundary of $\Omega$. For more, see \cite{MaulBanski1996}. The left shift, $S$, defined by $S(x_1,\cdots,x_n)=(x_2,\cdots,x_n)$, corresponds to applying $F$. Note that
$$
[X]=\bset{w\in\partial\Omega, \bb{\pi\inv(w)}_k=x_k, 1\le k\le \abs X}=f_{X|_{n-1}}\bb{\partial_{x_n}}.
$$
This makes sense since $F|_{U_j}=f_j\inv$ and when inverting the composition of functions, the order of composition is reversed.

\subsection{Refined Carleson's  Estimates}
The next lemma is a refinement of Carleson's estimate on the multiplicativity of harmonic measures of Jordan Repellers. It is a quantified version that will allow us to propagate the 'good disks' used to define $N(\delta,\alpha,\gamma;\eta)$ into smaller scales with a uniform error.

\begin{defn}
We say a Jordan Repeller \emph{expands at rate at least $D>1$} if for every mapping in the system $\abs{F'_j}\ge D$. Equivalently, for every disk of radius $r$, $B$, in $U_j$, $\diam(f_j(B))\le \frac rD$.
\end{defn}

We shall describe some notation that will be used in the next two lemmas; 

Let $F$ be a Jordan Repeller expanding at the rate of at least $D>1$, and let $Q_j:=U_j\cap \Omega_F$,  $\bset{U_j}$ are the neighbourhoods where $F|_{U_j}$ are invertible. Assume that the extremal distance between $\partial_j$ and $\partial Q_j\cap\Omega_F$ in $Q_j$ is uniformly bounded from above and below by $\frac1M$ and $M$ respectively.

Let $X=(x_1,\cdots,x_n)$ be a cylinder and denote by $U_X=f_{x_1}\circ f_{x_2}\circ\cdots\circ f_{x_{n-1}} U_{x_n}$, by $Q_X=f_{x_1}\circ f_{x_2}\circ\cdots\circ f_{x_{n-1}} Q_{x_n}$, and by $\lambda_X=f_{x_1}\circ f_{x_2}\circ\cdots\circ f_{x_{n-1}}\lambda_{x_n}$, where $\lambda_{j}\in Q_{j}$ are chosen so that $\lambda_j=\phi\bb{\bb{1-\omega\bb{\partial_j}}\zeta_j}$, and $\zeta_j\in \phi\inv\bb{\partial_j}$ is the center of the arc $\partial_j$. 
\begin{lem}[Refined Carleson's estimate]\label{lem:Carleson}
There exists a constant $C=C(M)$, so that
$$
\abs{\frac{\omega\bb{XYZ}}{\omega\bb{XY}}\cdot\frac{\omega(Y)}{\omega\bb{YZ}}-1}\le C\cdot \bb{\frac1D}^{\abs Y-1}.
$$
\end{lem}
We rely on Makarov's proof in \cite[p.52-53]{maksurvey}. We will need the following claim, which is an application of Koebe's distortion theorem.
\begin{claim}\label{clm:derivatives}
Let $\Omega\subset\C$ be a domain, and let $K\subset\Omega$ be a connected compact set of diameter 1 such that the modulus of the doubly connected domain $\Omega\setminus K$ is bounded from above and below by $M$ and $\frac1M$ respectively. Let $h:\Omega\rightarrow \D$ be a conformal map satisfying $0\in h(K)$. Then, for all $z\in K$, $\abs{\left.h'\right|_K}\asymp 1$ and $\abs{\left.h''\right|_K}\asymp 1$ with constants which depend on $M$.
\end{claim}

\begin{proof}[The proof of Lemma \ref{lem:Carleson}]
Let $g:\Omega_F\rightarrow R\D$ be a conformal map mapping $z_0$ to the origin for some $R$ large enough so that for every $j$, $\lambda_1(g(\partial_j))\ge 1$. For every $j$ we denote by $G_j=g(Q_j)$, $\alpha_j=g(\partial_j)$, and $\sigma_j=g(\partial Q_j\cap\Omega_F)$. Let $h_j:G_j\rightarrow\D^+$ be a conformal map which maps $\sigma_j$ to $\T^+$ and the centre of $\alpha_j$ to the origin. We then define $\nu_j:=g\inv \bb{h_j\inv\bb{\frac i2}}\in Q_j$, and denote by $\tilde\alpha_j=h_j(\alpha_j)$. See Figure \ref{fig:Carleson}.

\begin{figure}[htb]
\centering
\includegraphics[scale=0.8]{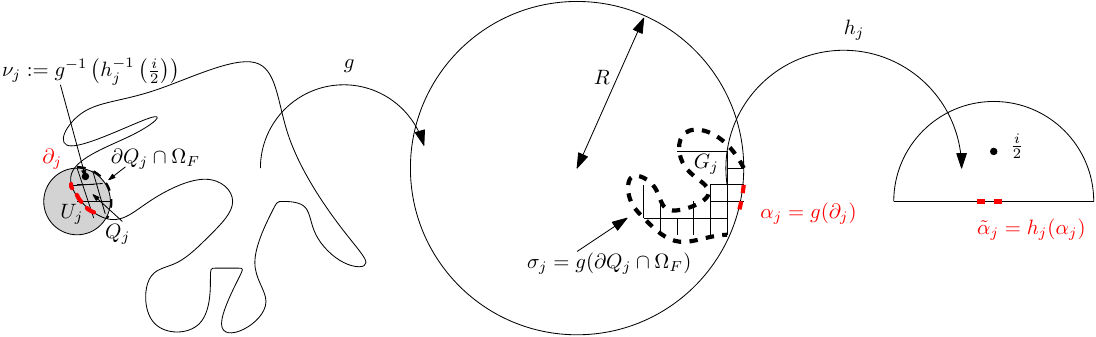}
\caption{This figure sketches the different sets and different conformal maps described above. Note that the arc $\tilde\alpha$ lies on the real line, where the harmonic measure of the half disk is linear.}
\label{fig:Carleson}
\end{figure}

We will show that
\begin{equation}\label{eq:ratio}
\abs{\frac{\omega(XYZ)}{\omega(XY)}\cdot\bb{\frac{\omega\bb{\nu_{Xy_1},\partial_{XYZ};Q_{Xy_1}}}{\omega\bb{\nu_{Xy_1},\partial_{XY};Q_{Xy_1}}}}\inv-1}\dless M \bb{\frac1D}^{\abs Y-1},
\end{equation}
i.e., the constant depends on $M$ alone. The same holds for $\frac{\omega(YZ)}{\omega(Y)}$ and $\frac{\omega\bb{\nu_{y_1},\partial_{YZ};Q_{y_1}}}{\omega\bb{\nu_{y_1},\partial_{Y};Q_{y_1}}}$ and as harmonic measure is conformal invariant, $\omega\bb{\nu_{Xy_1},\partial_{XYZ};Q_{Xy_1}}=\omega\bb{\nu_{y_1},\partial_{YZ};Q_{y_1}}$ thus concluding the proof.

Let $G=g(Q_{Xy_1})$ and define $h^X:G\rightarrow \D^+$ by $h^X(z):=h_{y_1}\bb{g\circ F_X\inv\circ g\inv}=h_{y_1}\bb{g\circ f_X\circ g\inv}$. Note that by definition,
$$
h^X\circ g(\nu_{Xy_1})=\frac i2,\; h^X\circ g(\partial_{Xy_1})=\tilde\alpha_{y_1}, \text{ and } h^X\bb{\sigma_{Xy_1}}=\T^+.
$$

Denote by $\alpha:=g(\partial_{Xy_1}),\; \beta:=g(\partial_{XY}),\;\gamma:=g(\partial_{XYZ})$ and by $\tilde\alpha,\tilde\beta$, and $\tilde\gamma$ the images of $\alpha,\beta$, and $\gamma$ under $h^X$.

Since harmonic measure is conformal invariant,
$$
\frac{\omega(XYZ)}{\omega(XY)}=\frac{\omega(0,\gamma;R\D)}{\omega(0,\beta;R\D)}=\frac{\lambda_1\bb{\gamma}}{\lambda_1\bb{\beta}}\;,\;\;\text{ and }\;\;\;\;\frac{\omega\bb{\nu_{Xy_1},\partial_{XYZ};Q_{Xy_1}}}{\omega\bb{\nu_{Xy_1},\partial_{XY};Q_{Xy_1}}}=\frac{\omega\bb{\frac i2,\tilde\gamma;\D^+}}{\omega\bb{\frac i2,\tilde\beta;\D^+}}\dsim M\frac{\lambda_1(\tilde\beta)}{\lambda_1(\tilde\gamma)},
$$
since $\tilde\beta,\tilde\gamma\subset\tilde\alpha$ and $\dist(\tilde\alpha,\partial\D)\dmore M 1$, as we assume the extremal distance between $\partial_j$ and $\partial Q_j\cap\Omega_F$ is in $\bb{\frac1M,M}$, and extremal length is conformally invariant.

Then \eqref{eq:ratio} is equivalent to showing that
\begin{equation}\label{eq:new_ratio}
\abs{\frac{\lambda_1\bb{\gamma}}{\lambda_1\bb{\beta}}\bb{\frac{\lambda_1\bb{\tilde\gamma}}{\lambda_1\bb{\tilde\beta}}}\inv-1}\dless M\bb{\frac1D}^{\abs Y-1}.
\end{equation}

Let $\hat G$ denote the symmetrization of $G$ across $\alpha$, and consider $\hat h^X$ as a map of the symmetrization, $\hat h^X:\hat G\rightarrow\D$. Because the extremal distance between $\partial_{y_1}$ and $ \partial Q_{y_1}\cap\Omega$ is assumed to be bounded between $\frac 1M$ and $M$, the domain $\hat G$ and the compact set $\alpha$ satisfy the requirement of Claim \ref{clm:derivatives},
$$
\abs{\left.\bb{h^X}'\right|_\alpha}\dsim M1,\;\;\abs{\left.\bb{h^X}''\right|_\alpha}\dsim M 1.
$$
Fix $\zeta_0\in\beta$, then since $\tilde\beta$ is the image of $\beta$ under the map $h^X$,
\begin{align*}
\abs{\lambda_1\bb{\tilde\beta}-\lambda_1\bb{\beta}\abs{\bb{h^X}'(\zeta_0)}}&=\abs{\underset\beta\int\abs{\bb{h^X}'(\zeta)}d\abs{\zeta}-\lambda_1\bb{\beta}\abs{\bb{h^X}'(\zeta_0)}}\\
&\le \underset\beta\int\abs{\bb{h^X}'(\zeta)-\bb{h^X}'(\zeta_0)}d\abs{\zeta}=\underset\beta\int\abs{\bb{h^X}''(\xi_\zeta)}\abs{\zeta-\zeta_0}d\abs{\zeta}\dless M\lambda_1\bb{\beta}^2\\
&\Rightarrow\abs{\lambda_1\bb{\beta}\inv\lambda_1\bb{\tilde\beta}-\abs{\bb{h^X}'(\zeta_0)}}\dless M \lambda_1\bb{\beta}.
\end{align*}

The same argument done with the curve $\gamma$ shows that
$$
\abs{\lambda_1\bb{\gamma}\inv\lambda_1\bb{\tilde\gamma}-\abs{\bb{h^X}'(\zeta_0)}}\dless M\lambda_1\bb{\gamma}.
$$
This implies that there exists a constant $C=C(M)$ so that
$$
\abs{\frac{\lambda_1\bb{\gamma}}{\lambda_1\bb{\beta}}\bb{\frac{\lambda_1\bb{\tilde\gamma}}{\lambda_1\bb{\tilde\beta}}}\inv-1}\le \abs{\frac{\abs{\bb{h^X}'(\zeta_0)}+C\lambda_1\bb{\beta}}{\abs{\bb{h^X}'(\zeta_0)}-C\lambda_1\bb{\gamma}}-1}=\frac{C\bb{\lambda_1\bb{\beta}+\lambda_1\bb{\gamma}}}{\abs{\bb{h^X}'(\zeta_0)}-C\lambda_1\bb{\gamma}}\dless M\lambda_1\bb{\beta}.
$$
To conclude the proof of \eqref{eq:new_ratio} it is left to show that $\lambda_1(\beta)\dless M\bb{\frac1D}^{\abs Y-1}$.

We note that since the extremal distance between $\frac i2$ and $\partial\D^+$ is some constant (for definition see \cite[p.144 bottom]{gm}), then extremal distance between $\partial Q_{y_1}$ and $\nu_{y_1}$ is uniformly bounded from above and bellow by  numerical constants as a conformal image of the half disk. This implies that for every arc in the boundary of $\partial Q_{y_1}$, $\omega\bb{\nu_{y_1},A;Q_{y_1}}\sim\frac{\lambda_1(A)}{\diam(Q_{y_1})}$, where the constant is a numerical constant. Overall
\begin{align*}
\lambda_1\bb{\beta}\sim\omega\bb{\frac i2,\tilde\beta;\D^+}=\omega\bb{\nu_{Xy_1},\partial_{XY};Q_{Xy_1}}&=\omega\bb{\nu_{y_1},\partial_{Y};Q_{y_1}}\dsim M\frac{\lambda_1\bb{\partial_Y}}{\diam(Q_{y_1})}\\
&=\frac{\lambda_1\bb{f_{y_n}\circ\cdots\circ f_{y_2} \partial_{y_1}}}{\diam(Q_{y_1})}\dless M \bb{\frac1D}^{\abs Y-1}.
\end{align*}

To conclude the proof, it is left to show that
$$
\abs{\frac{\lambda_1\bb{\tilde\gamma}}{\lambda_1\bb{\tilde\beta}}:\frac{\omega\bb{\frac i2,\tilde\gamma;\D^+}}{\omega\bb{\frac i2,\tilde\beta;\D^+}}-1}\lesssim\bb{\frac1D}^{\abs Y-1}.
$$
Note that while we know that $\lambda_1(\tilde\beta)\sim\omega\bb{\frac i2,\tilde\beta;\D^+}$, this is not fine enough for the estimate we need.

Let $\phi:\D^+\rightarrow\D$ be a conformal map, mapping $\frac i2$ to the origin, and $\T^+$ to itself. Then for every arc $A\subset\partial\D^+$,
$$
\omega\bb{\frac i2,A;\D^+}=\omega(0,\phi(A);\D)=\lambda_1(\phi(A)).
$$
Now, $\phi$ is a fixed M\"obius map, and therefore its second derivative is uniformly bounded as a function of the distance of $\tilde\alpha$ from $\pm 1$, which is equal (up to a uniform constant) to the extremal distance between $\partial_{y_1}$ and $\partial Q_{y_1}\cap\Omega$, in $Q_{y_1}$, in other words, it depends on $M$. Fix $\tilde\zeta_0\in\tilde\gamma$, then , as before, for every arc $A$ containing $\tilde\gamma$,
\begin{align*}
\abs{\lambda_1(\phi(A))-\lambda_1(A)\abs{\phi'(\tilde\zeta_0)}}&\le\int_A\abs{\abs{\phi'(\zeta)}-\abs{\phi'(\tilde\zeta_0)}}d\abs\zeta\le \int_A\abs{\phi'(\zeta)-\phi'(\tilde\zeta_0)}d\abs\zeta\\
&= \int_A\abs{\phi''(\xi_\zeta)}\abs{\zeta-\tilde\zeta_0}d\abs\zeta\lesssim \lambda_1(A)^2,
\end{align*}
where the constant only depends on the the bounds we had for the second derivative, which in turn depends on $M$. Overall,
\begin{align*}
\abs{\frac{\omega\bb{\frac i2,\tilde\gamma;\D^+}}{\omega\bb{\frac i2,\tilde\beta;\D^+}}\cdot\bb{\frac{\lambda_1\bb{\tilde\gamma}}{\lambda_1\bb{\tilde\beta}}}\inv-1}&=\abs{\frac{\lambda_1(\phi(\tilde\gamma))}{\lambda_1(\phi(\tilde\beta))}\cdot\bb{\frac{\lambda_1\bb{\tilde\gamma}}{\lambda_1\bb{\tilde\beta}}}\inv-1}\le \frac{\bb{\lambda_1(\tilde\gamma)\abs{\phi'(\tilde\zeta_0)}+\lambda_1(\tilde\gamma)^2}\lambda_1\bb{\tilde\gamma}\inv}{\bb{\lambda_1(\tilde\beta)\abs{\phi'(\tilde\zeta_0)}+\lambda_1(\tilde\beta)^2}\lambda_1\bb{\tilde\beta}\inv}-1\\
&\le\frac{\abs{\phi''(\tilde\zeta_0)}\inv\bb{\lambda_1(\tilde\gamma)+\lambda_1(\tilde\beta)}}{1+\abs{\phi''(\tilde\zeta_0)}\inv\lambda_1(\tilde\beta)}\dless M\lambda_1(\tilde\beta)\dless M\lambda_1(\beta)\lesssim\bb{\frac1D}^{\abs Y-1},
\end{align*}
concluding the proof.
\end{proof}

\begin{rmk}
For rotations the story is slightly different. The reason is that Jordan Repellers might be defined only near the boundary, but the rotation passes through the entire domain. This is why the very first step, which brings us from the marked point, $z_0$, to the domain where the Jordan Repeller, $F$, is defined, could have a very large rotation, which does not repeat later in the iterations. Another reason this is different from harmonic measure is that the harmonic measure is invariant under conformal maps, while rotation is not.
\end{rmk}
\begin{rmk}\label{rmk:rotation_curve}
In the context of Jordan Repellers, it is convenient to talk about the rotation of a finite word. For a finite word, $X$, we denote by $\rot(X)$ the rotation of the crosscut supported on $\phi\inv\bb{\varphi[X]}$. Note that for quasidisks, for every arc in $\partial\D$, $\phi(C)$ is a crosscut, and so a crosscut supported on $\phi\inv\bb{ [X]}$ always exists.
\end{rmk}
\begin{lem}[Carleson's type estimate for rotations]\label{lem:Carleson_rotations}
Let $F$ be a Jordan Repeller, and fix $\abs X=n, \abs Y=k$, $n,k\ge 2$ be two finite words. Then $\frac{\rot(XY)}{\rot\bb{\partial_{y_{k}}}}\sim\frac{\rot(X)}{\rot\bb{\partial_{x_{n}}}}\cdot\frac{\rot(Y)}{\rot\bb{\partial_{y_{k}}}}$ and the constants are numerical.
\end{lem}
\begin{proof}
For every word $X=(x_1,x_2,\cdots)$ we let $\lambda_X\in Q_X$ be the point satisfying that $F^{\circ \bb{\abs X-1}}(\lambda_X)=\lambda_{x_n}$, and set $z_X:=\phi\inv\bb{\lambda_X}$. Following Lemma \ref{lem:rotation_differentiation}, if $1-\abs{z_X}\sim\omega(\partial_X)$, then $\abs{\phi'(z_X)^{-i}}\sim\rot(X).$ We shall first show that $1-\abs{z_X}\sim\omega(\partial_X)$ by induction. For $X=x_1$, this is true by definition. Assume this is the case for all words of length at most $(n-1)$ and let $X$ be a word of length $n$. Then $X=x_1Y$, where $Y$ is a word of length $(n-1)$, and $F(\lambda_X)=\lambda_Y$.

Define the map $B:\bb{\bunion j 1 N U_j}\cap\Omega\rightarrow\C$ by $B(z)=\phi\inv\circ F\circ \phi$, then $B^{\circ n}(z)=\phi\inv\circ F^{\circ n}\circ \phi$ and, in particular, 
$$
B(z_X)=\phi\inv\circ F\circ \phi (z_X)=\phi\inv\circ F(\lambda_X)=\phi\inv\circ \lambda_Y=z_Y.
$$
In addition, since $\phi$ and $F$ are both conformal, so is $B$, which means it preserves extremal length. Recall that $U_Y$ is a polar cube centered at $z_Y$. This implies that $U_X=B\inv (U_Y)$ is a polar cube centered at $z_X=B\inv(z_Y)$.

Next, following Lemma \ref{lem:rotation_differentiation},  for every finite word, $X$, 
$$
\rot(X)=\exp(\log\rot(X))= \exp\bb{\Im\bb{\log \phi'(z_X)}+C}=C_1\cdot \abs{\bb{\phi'(z_X)}^{-i}}.
$$
Recall that $S$ denotes the left shift. Then for every finite word, $X$, with $n=\abs X$,
$$
\phi'(z_X)=\phi'(z_{S^{\circ 0}(X)})=\prodit j 0 {n-1}\frac{\phi'(z_{S^{\circ j}(X)})}{\phi'(z_{S^{\circ {(j+1)}}(X)})}\cdot \phi'(z_{S^{\circ {n-1}}(X)})=\phi'(z_{x_n})\prodit j 0 {n-1}\frac{\phi'(z_{S^{\circ j}(X)})}{\phi'(z_{S^{\circ {(j+1)}}(X)})}.
$$

We see that if $\abs X=n, \abs Y=k$, then
\begin{align*}
\phi'\bb{z_{XY}}=\phi'\bb{z_{y_{k}}}\cdot\prodit j 0 {n+k-1}\frac{\phi'\bb{z_{S^{\circ j}(XY)}}}{\phi'\bb{z_{S^{\circ \bb{j+1}}(XY)}}}=\phi'\bb{z_{y_{k}}}\cdot\bb{\prodit j 0 {n-1}\frac{\phi'\bb{z_{S^{\circ j}(XY)}}}{\phi'\bb{z_{S^{\circ \bb{j+1}}(XY)}}}}\cdot\bb{\prodit j {n} {n+k-1}\frac{\phi'\bb{z_{S^{\circ j}(XY)}}}{\phi'\bb{z_{S^{\circ \bb{j+1}}(XY)}}}}\\
=\bb{\prodit j 0 {n-1}\frac{\phi'\bb{z_{S^{\circ j}(XY)}}}{\phi'\bb{z_{S^{\circ \bb{j+1}}(XY)}}}}\cdot\bb{ \phi'\bb{z_{y_{k}}}\prodit j 0 {k-1}\frac{\phi'\bb{z_{S^{\circ j}(S^{\circ \bb{n}}(XY))}}}{\phi'\bb{z_{S^{\circ \bb{j+1}}(S^{\circ \bb{n}}(XY))}}}}\\
=\frac1{\phi'\bb{z_{x_{n}}}}\bb{\phi'\bb{z_{x_{n}}}\prodit j 0 {n-1}\frac{\phi'\bb{z_{S^{\circ j}(XY)}}}{\phi'\bb{z_{S^{\circ \bb{j+1}}(XY)}}}}\cdot \bb{\phi'\bb{z_{y_{k}}}\cdot\prodit j 0 {k-1}\frac{\phi'\bb{z_{S^{\circ j}(Y)}}}{\phi'\bb{z_{S^{\circ \bb{j+1}}(Y)}}}}
\end{align*}
by modifying the index in the second sum. Using Lemma \ref{lem:rotation_differentiation} which relates the rotation with the derivative,
\begin{align*}
\frac{\rot(XY)}{\rot\bb{\partial_{y_{k}}}}&\sim\abs{\bb{\frac{\phi'\bb{z_{XY}}}{\phi'\bb{z_{y_{k}}}}}^{-i}}\\
&= \frac1{\abs{\bb{\phi'\bb{z_{x_{n}}}}^{-i}}}\abs{\bb{\phi'\bb{z_{x_{n}}}\prodit j 0 {n-1}\frac{\phi'\bb{z_{S^{\circ j}(XY)}}}{\phi'\bb{z_{S^{\circ \bb{j+1}}(XY)}}}}^{-i}}\cdot \frac1{\abs{\bb{\phi'\bb{z_{y_{k}}}}^{-i}}}\cdot\abs{\bb{\phi'\bb{z_{y_{k}}}\cdot\prodit j 0 {k-1}\frac{\phi'\bb{z_{S^{\circ j}(Y)}}}{\phi'\bb{z_{S^{\circ \bb{j+1}}(Y)}}}}^{-i}}\\
&=\frac{\abs{\bb{\phi'\bb{z_{X}}}^{-i}}}{\abs{\bb{\phi'\bb{z_{x_{n}}}}^{-i}}}\cdot\frac{\abs{\bb{\phi'\bb{z_{Y}}}^{-i}}}{\abs{\bb{\phi'\bb{z_{y_{k}}}}^{-i}}}\sim \frac{\rot(X)}{\rot\bb{\partial_{x_{n}}}}\cdot\frac{\rot(Y)}{\rot\bb{\partial_{y_{k}}}}
\end{align*}
with some numerical constant.
\end{proof}

\subsection{Spectra of words}
Let $\Sigma:=\bset{a_1,\cdots,a_M}$ be a finite alphabet describing a Jordan Repeller, $F$, and let $d:=\underset {a\in\Sigma}\min\; \diam([a]),\; D:=\underset {a\in\Sigma}\max\; \diam([a])$.

\begin{rmk}\label{rmk:renormalised}
We would like to present similar definitions for the dimension and the distortion spectra in the context of words. However, while the initial diameter of the repeller has no effect on the harmonic measure or the rotation, as these are scale-invariant, it has a grave impact on the multiplicativity property of the diameter. Looking at the case of a Carleson fractal, that is, when all maps are linear, the magnitude of the derivative of each map is $\frac{\diam([a_j])}{\diam(\partial)}$, which implies that the diameter of a word compared to the diameter of the letters composing it is
$$
\frac{\diam([a_1a_2])}{\diam([a_1])\cdot\diam([a_2])}=\frac{\diam([a_1])\cdot \frac{\diam([a_2])}{\diam(\partial)}}{\diam([a_1])\cdot\diam([a_2])}=\frac1{\diam(\partial)}.
$$
For a word composed of $n$ letters, the latter should be raised to a power of $n-1$. This demonstrates that some normalization is required for the analysis of the spectra, where the multiplicativity property of the harmonic measure and of the rotation should go hand in hand with the multiplicativity of the diameter.

A similar effect is obtained by the rotation and the harmonic measure of the repeller. We therefore define the spectra of the repeller using these normalised quantities.

We will abuse the notation of harmonic measure, rotation, and diameter from now on, to refer to their renormalised counterparts.
\end{rmk}
\paragraph{Definitions:}
Let $F$ be a Jordan Repeller described by a symbolic dynamical system with alphabet $\Sigma=\bset{a_1,\cdots,a_N}$. Following Remark \ref{rmk:renormalised}, we abuse the notation of diameter, harmonic measure, and rotation of words by defining 
\begin{equation}\label{eq:renormalised}
\diam(w)=\frac{\diam([w])}{\diam\partial}\quad;\quad\omega(w)=\frac{\omega(z_0,[w];F)}{\omega(z_0,\partial;F)}\quad;\quad \rot(w)=\frac{\rot\bb{z_0, [w];F}}{\rot\bb{z_0,\partial;F}},
\end{equation}
where since $ [w]$ is a curve, the definition of rotation of a finite word, relies on the definition of rotation of a crosscut (see Remark \ref{rmk:rotation_curve}). In addition, we denote by $\abs w$ the length of the word $w$.

For every $\delta>0$  we denote by $I^\delta$ the collection of all indices giving rise to finite words with diameter $\delta$:
$$
I^\delta:=\bset{(k_1,\cdots,k_N);\prodit j 1N  \bb{\diam\bb{a_j}}^{k_j}=\delta, N\in\N}.
$$
For example, if all the letters have the same diameter, then
$$
I^\delta=\bset{(k_1,\cdots,k_N);  \delta=\diam(a_1)^{\sumit j1 N k_j}}.
$$
However, in the most general case, this set could vary wildly depending on the diameters of the letters. For example, it could be the case that $\diam(a_1)\cdot \diam(a_2)=\diam(a_3)$ and therefore this collection will include words of different length.

Given a multi-index $\vect i=(k_1,\cdots,k_N)\in I^\delta$ we say $w=(w_1,\cdots,w_m)\in W^{\vect i}$ if for every $1\le j\le N$
$$
\#\bset{\nu, w_\nu=a_j}=k_j.
$$
\begin{defn}\label{defn:word_spectra}
We define the \emph{Minkowski Word Dimension Spectrum} by 
$$
f_\Omega^{\sigma\sigma'\;word}(\alpha,\gamma)=\Ulimit \eta 0 \Ulimitsup \delta 0 \frac{\log N_{word}^{\sigma,\sigma'}(\delta,\alpha,\gamma,\eta)}{\log\bb{\frac1\delta}},\quad \sigma,\sigma'\in\bset{+,-}
$$
where $N_{word}^{\sigma\sigma'}(\delta,\alpha,\gamma,\eta)$ is the maximal number of disjoint words $w\in \underset{\vect i\in I^\delta}\bigcup W^{\vect i}$,  satisfying \eqref{eq:N_measure} and \eqref{eq:N_rotation} (in their renormalised version appearing in \eqref{eq:renormalised}). 
\end{defn}

The next natural definition to extend is the Minkowski Distortion Word Spectrum. However, since $ [w]$ is a curve, it is clear that this quantity is dwarfed by the Minkowski Crosscut Spectrum. To see it is in fact the same, for every curve $\nu$ we define $\underline w_\nu$ to be the maximal word satisfying $ [\underline w_\nu]\subseteq\nu$ and $\overline w_\nu$ to be the minimal word satisfying $ [\overline w_\nu]\subseteq\nu$. Since the number of letters distinguishing between $\overline w_\nu$ and $\underline w_\nu$ is bounded by four, and following Carleson's refinement, we see that
$$
\frac{\omega([\overline w_\nu])}{\omega([\underline w_\nu])}\in\left[1,C\right)\quad, \quad\frac{\rot([\overline w_\nu])}{\rot([\underline w_\nu])}\in\bb{C\inv,C}\quad,\quad \frac{\diam([\overline w_\nu])}{\diam([\underline w_\nu])}\in\left[1,C\right).
$$
Then following the same argument as the one presented in the proof of Lemma \ref{lem:curves_bnd_d} part \ref{item:quasidisk_curve}, we see that the number of such words corresponds to the measure of $L_{a',b'}(r')$ and the later (in the limit) is equal to the Minkowski Crosscut Spectrum, following Lemma \ref{lem:curves_bnd_d} part \ref{item:quasidisk_curve} as Jordan Repellers are quasidisks.

\paragraph{Consistency:}
\begin{lem}\label{lem:consistent_defs} Let $F$ be a finite Jordan Repeller. Then
\begin{equation}\label{eq:equivalent_f}
f_F^{\sigma,\sigma'\;word}(\alpha,\gamma)=f_F^{\sigma,\sigma'}(\alpha,\gamma)\quad,\quad \sigma,\sigma'\in\bset{+,-}.
\end{equation}
where $f_F^{\sigma,\sigma'}(\alpha,\gamma)$ is the Minkowski Dimension Mixed Spectrum with disks centered at $\partial$.
\end{lem}
We begin by proving a more general statement that will be used to show consistency.
\begin{prop}\label{prop:good-vs-good}
Let $F$ be a Jordan Repeller expanding at a rate at least 2.
\begin{enumerate}[label=(\roman*)]
\item\label{itm:1} For every finite word, $w$, let $B=B_w$ be the smallest disk containing $ ([w])$. Then
\begin{align*}
\frac{\diam(B)}{\diam( [w])}\in\bb{1,2}\quad\quad,\quad\quad \omega([w])\sim\omega(B)\quad\quad;\quad\quad\rot([w])\sim \rot(B).
\end{align*}
\item\label{itm:2} For every disk, $B$, and $\eta>0$ there exists a word, $w$, and $k\in\N$ satisfying
\begin{align*}
&r^{1+\eta}\le \diam( [w]_0^k)\le r\\
&\omega(B)^{1+2\eta}\le \omega( [w]_0^k)\le \omega(B)^{1-\eta}\quad\quad;\quad\quad\rot(B)^{1+2\eta}\le \rot( [w]_0^k)\le \rot(B)^{1-\eta}.
\end{align*}
\end{enumerate}
\end{prop}
\begin{proof}
To see the first statement, we use Lemma \ref{lem:doubling_measure_quasidisk} and the fact that there are only finitely many components of $2 B$ intersecting $B$ (by Claim \ref{clm:quasi_finite_curves}). Then $\omega([w])\sim\omega(B)$. Since the diameter of the curve $ [w]$ and the disk are comparable and these sets intersect, $\rot([w])\sim \rot(B)$, following Lemma \ref{lem:rotation_err_in}, concluding the proof of \ref{itm:1}.

To see \ref{itm:2}, fix a disk, $B=B(z,r)$, and $\eta>0$. Because Jordan Repellers are quasidisks, their harmonic measure is doubling. Therefore, we may assume without loss of generality that there exists an infinite word so that $ [w]=z$, is the centre of the disk, $B$ (otherwise, because this correspondence is defined almost surely, we can shift the disk slightly and change the harmonic measure by at most a constant). Let
$$
k_{min}:=\min\bset{\nu\in\N, \diam( [w]_0^\nu)\le r}.
$$
On one hand,
$$
r\ge \diam( [w]_0^k)\ge \diam( [w]_0^{k-1})\cdot\underset{a_j\in\Sigma}\min \diam(a_j)>\frac{\diam( [w]_0^{k-1})}2 >\frac r2\ge r^{1+\eta}
$$
as long as $r$ is small enough (depending on $\eta$). On the other hand, because $ [w]\in B\cap  [w]_0^k$, then
$$
2B=B( [w],2r)\supseteq B( [w],\diam( [w]_0^{k-1}))\supseteq  [w]_0^k.
$$
Next, the harmonic measure is doubling, implying that
$$
\omega( [w]_0^k)\le \omega(2B)\sim\omega(B).
$$
Similarly, 
$$
  [w]_0^k\supseteq B\bb{ [w],\frac12\cdot \diam( [w]_0^{k-1})}\supseteq B\bb{ [w],\frac12 r^{1+\eta}},
$$
which implies that
$$
\omega(  [w]_0^k)\gtrsim\omega(B)^{1+\eta}.
$$
As before, since the diameter of the curve $ [w]$ and the disk are comparable and these sets intersect, $\rot([w])\sim \rot(B)$, following Lemma \ref{lem:rotation_err_in}.
\end{proof}

\begin{proof}[Proof of Lemma \ref{lem:consistent_defs}]

To conclude equality \eqref{eq:equivalent_f}, it is left to note that every sequence $\delta_\nu\searrow 0$ corresponds to a sequence of configurations $(k_1^\nu,\cdots k_N^\nu)\in I^{\delta_\nu}$ satisfying 1, 2, and 3 giving rise to a collection of disks counted in the definition of the Minsowski Dimension Spectrum, i.e., their harmonic measure and its rotation satisfy \eqref{eq:N_measure} and \eqref{eq:N_rotation} respectively, and vise-verse; every collection of disks, $\bset{B_j}$, counted in the Minsowski Dimension Spectrum, give rise to configurations $(k_1^\nu(j),\cdots k_N^\nu(j))\in I^{\delta_\nu}$ satisfying 1, 2, and 3. While disjoint disks correspond to distinct words, distinct words could potentially give rise to intersecting disks. However, by excluding a linear portion of the words, one could guarantee that the disks arising from them are indeed disjoint, concluding the proof.
\end{proof}

\subsection{Propagation:}
For the reader's convenience we repeat the definition of a propagating function here:
\begin{defn*}
A function $ :(0,1)^2\rightarrow\R$ \emph{propagates} if there exists a constant $C$, so that for every $n\in\N$
$$
\frac{\log \varphi\bb{\delta^{n},C\cdot \eta}}{\log\bb{\frac1{\delta^{n}}}}\ge \frac{\log \varphi(\delta,\eta)}{\log\bb{\frac1\delta}}
$$
for every $\delta$ small enough (which may depend on $\eta$).
\end{defn*}
The first observation is that for Jordan Repellers, the functions $\varphi(\delta,\eta):=N_{word}^{\sigma,\sigma'}\bb{\delta,\alpha,\gamma,\eta},\;\sigma,\sigma'\in\bset{+,-}$ and $\varphi((1-r),\abs{a'-a})=\#\Gamma(a',b,r)$ propagate, making this property interesting.
\begin{lem}\label{lem:propagations}
Let $F$ be a Jordan Repeller.
\begin{enumerate}
\item \label{itm:prop_N} The function $(\delta,\eta)\mapsto N_{word}^{\sigma,\sigma'}(\delta,\alpha,\gamma,\eta),\;\sigma,\sigma'\in\bset{+,-}$ propagates.
\item \label{itm:prop_curves}The function $(\delta,\eta)\mapsto \#\Gamma(a-\eta,b-\eta,1-\delta)$ propagates.
\end{enumerate}
\end{lem}
\begin{rmk}~
\begin{enumerate}[label=(\arabic*)]
\item In fact, following the proof, we will see that the constant that appears in the definition of propagation is the same constant given by Carleson's estimates, Lemmas \ref{lem:Carleson} and \ref{lem:Carleson_rotations}.
\item If there exists a collection of words, $\mathcal M$, satisfying
\begin{enumerate}[label=$\bullet$]
\item For every $w,w'\in\mathcal M$, $dist( [w], [w'])>\delta$
\item $\diam( [w])<\delta$
\item $ [w]$ satisfies \eqref{eq:N_measure}, and \eqref{eq:N_rotation}
\end{enumerate}
then $N_{word}^{\sigma,\sigma'}(\delta,\alpha,\gamma,\eta)\gtrsim\#\mathcal M$, and since distances are rescaled accordingly,
$$
\frac{\log N_{word}^{\sigma,\sigma'}(\delta^n,\alpha,\gamma,C\cdot \eta)}{\log\bb{\frac1{\delta^{n}}}}\ge \frac{\log \#\mathcal M}{\log\bb{\frac1\delta}}-o\bb{1}.
$$
\end{enumerate}
\end{rmk}
We will use the following three observations:

\paragraph*{Observation 1:} If $w_1,\cdots w_n\in\underset{\vect i\in I^\delta}\bigcup W^{\vect i}$, then for every permutation $\tau$, 
$$
P(w_1,\cdots,w_n,\tau):=w_{\tau(1)}w_{\tau(2)}\cdots w_{\tau(n)}\in \underset{\vect {i'}\in I^{\delta^n}}\bigcup W^{\vect {i'}}.
$$
Assume that for every word $w_j$, let $\vect i_j\in I^\delta$ be so that $w_j\in W^{\vect i_j}$, then if $\vect i_j=(k_1^j,\cdots,k_N^j)$ then
$$
\diam\bb{P(w_1,\cdots,w_n,\tau)}=\prodit j 1 n \prodit \ell 1 N \bb{a_\ell}^{k_\ell^j}=\prodit j 1 n\delta=\delta^n
$$
or, in other words
$$
\bset{\vect i_j}_{j=1}^n\subset I^\delta\Rightarrow \sumit j 1 n \vect i_j\in I^{\delta^n}.
$$
\paragraph*{Observation 2:} If $w_1,w_2,\cdots, w_n\in\underset{\vect i\in I^\delta}\bigcup W^{\vect i}$, and $\tau$ is a permutation over $n$ elements, then by Carleson's refinement lemma, there exists $A>1$, which depends on a uniform upper bound on the extremal distance between $\partial_j$ and $\partial Q_j\cap\Omega_F$ alone, so that
$$
\frac{\omega\bb{z_0,\sbb{w_{\tau(1)}w_{\tau(2)}\cdots w_{\tau(n)}};F}}{\prodit j1 n\omega\bb{z_0,\sbb{w_j};F}}\in\bb{\frac1{A^n},A^n}\quad\quad;\quad\quad \frac{\rot\bb{z_0,\sbb{w_{\tau(1)}w_{\tau(2)}\cdots w_{\tau(n)}};F}}{\prodit j1 n\rot\bb{z_0,\sbb{w_j};F}}\in\bb{\frac1{A^n},A^n}.
$$
\paragraph*{Observation 3:} For every word $w=(a_1^w,a_2^w,\cdots,a_{\ell(w)}^w)\in W^{\vect i}$ we cannot change the length of $w$ and keep its diameter. In other words, if $S$ is the left shift operator, then
$$
Sw=(a_2^w,\cdots,a_{\ell(w)}^w)\nin \underset{\vect {i'}\in I^\delta}\bigcup W^{\vect {i'}}\quad\quad;\quad\quad\forall a\in\Sigma,\quad aw=(a,a_1^w,a_2^w,\cdots,a_{\ell(w)}^w)\nin \underset{\vect {i'}\in I^\delta}\bigcup W^{\vect {i'}}.
$$
This implies that if $w_1,w_2,\cdots, w_n\in\underset{\vect i\in I^\delta}\bigcup W^{\vect i}$, and $\tau$ is a permutation over $n$ elements, then there exists a \underline{unique} string 
$$
P(w_1,\cdots,w_n,\tau)\in \underset{\vect {i'}\in I^{\delta^n}}\bigcup W^{\vect {i'}},
$$
where the uniqueness is up-to repetition of words, i.e., the case where $w_k=w_j$ for some $j\neq k$.

\begin{proof}[The proof of Lemma \ref{lem:propagations}]~\\

{\bf The proof of \ref{itm:prop_N}:} Fix $\delta, \eta,\alpha,\gamma$ and let $\mathcal A(\delta, \eta,\alpha,\gamma)$ denote a maximal collection of words satisfying \eqref{eq:N_measure} and \eqref{eq:N_rotation}. For every $w_1,w_2,\cdots,w_n\in \underset{\vect i\in I^\delta}\bigcup W^{\vect i}$ and every permutation over $n$ elements, $\tau$, the word $P(w_1,\cdots,w_n,\tau)$ satisfies
\begin{enumerate}
\item Following Observation 1, $\diam\bb{\sbb{P(w_1,\cdots,w_n,\tau)}}=\delta^n$.
\item Following Observation 2,
\begin{align*}
&\omega\bb{\sbb{P(w_1,\cdots,w_n,\tau)}}		\begin{cases}
												\ge\frac1{A^n}\prodit j 1 n \omega\bb{\sbb{w_j}}\ge \frac1{A^n}\bb{\delta^{\alpha+\eta}}^n=\bb{\delta^n}^{\alpha+2\eta}\\
												\le A^n\prodit j 1 n \omega\bb{\sbb{w_j}}\le A^n\bb{\delta^{\alpha-\eta}}^n=\bb{\delta^n}^{\alpha-2\eta}
											\end{cases}\\
&\rot\bb{\sbb{P(w_1,\cdots,w_n,\tau)}}		\begin{cases}
												\ge\frac1{A^n}\prodit j 1 n \rot\bb{\sbb{w_j}}\ge \frac1{A^n}\bb{\delta^{\gamma+\eta}}^n=\bb{\delta^n}^{\gamma+2\eta}\\
												\le A^n\prodit j 1 n \rot\bb{\sbb{w_j}}\le A^n\bb{\delta^{\gamma-\eta}}^n=\bb{\delta^n}^{\gamma-2\eta}
											\end{cases}											
\end{align*}
if $\delta$ is small enough, depending on $A$ and $\eta$.
\item Following Observation 3, if $w_1,\cdots, w_m\in \mathcal A(\delta, \eta,\alpha,\gamma)$ and we use the $j$th word $t_j$ times with $\sumit j 1 m t_j=n$, then there are 
$$
{n\choose t_1,t_2,\cdots t_m}=\frac{n!}{t_1! t_2!\cdots t_m!}=\frac{\#\bset{\text{permutations over n elements}}}{\#\bset{\text{orders of copies of the same word}}}
$$
different distinct words that belong to $\mathcal A(\delta^n, 2\eta,\alpha,\gamma)$.
\end{enumerate}
We conclude that for every $n$,
\begin{align*}
N_{word}^{\sigma,\sigma'}(\delta^n,\alpha,\gamma,2\eta)&\ge \sumit m 1 {N_{word}^{\sigma,\sigma'}(\delta,\alpha,\gamma,\eta)}\underset{(w_1,\cdots,w_m)\in \mathcal A(\delta, \eta,\alpha,\gamma)}\sum\underset{(t_1,\cdots,t_m)\atop\sumit j 1 m t_j=n}\sum{n\choose t_1,t_2,\cdots t_m}\\
&=\sumit m 1 {N_{word}^{\sigma,\sigma'}(\delta,\alpha,\gamma,\eta)}{N_{word}^{\sigma,\sigma'}(\delta,\alpha,\gamma,\eta)\choose m}\cdot m^n\ge \bb{N_{word}^{\sigma,\sigma'}(\delta,\alpha,\gamma,\eta)}^n
\end{align*}
using a theorem on multinomial coefficients. We conclude that

$$
\frac{\log \bb{N_{word}^{\sigma,\sigma'}(\delta^n,\alpha,\gamma,2\eta)}}{\log\bb{\frac1{\delta^{n}}}}\ge \frac{\log \bb{\bb{N_{word}^{\sigma,\sigma'}(\delta,\alpha,\gamma,\eta)}^n}}{n\log\bb{\frac1{\delta}}}=\frac{\log \bb{N_{word}^{\sigma,\sigma'}(\delta,\alpha,\gamma,\eta)}}{\log\bb{\frac1{\delta}}},
$$
i.e., the function $N_{word}^{\sigma,\sigma'}(\delta,\alpha,\gamma,\eta)$ propagates.

The proof of  \ref{itm:prop_curves} is identical, and therefore we leave it for the reader.
\end{proof} 

\begin{cor}\label{cor:propagation}
If the function $\varphi=\varphi_\Omega$ propagates, then for every $\eps>0$ there exists $\delta_0>0$ and $\eta_0>0$ so that
$$
\limit\eta0\limit\delta0\frac{\log \varphi(\delta,\eta)}{\log\bb{\frac1\delta}}\ge \frac{\log \varphi(\delta_0,\eta_0)}{\log\bb{\frac1{\delta_0}}}-\eps.
$$
In particular, one can approximate the spectra of Jordan Repellers by a finite scale spectra.
\end{cor}